\documentclass[MSNbibl,number,citesort,dvips]{arxbj}
\usepackage{upgreek}
\usepackage{dcolumn}
\usepackage[left]{vmkcol}

\aid{0}
\volume{20}
\issue{1}
\pubyear{2014}
\firstpage{141}
\lastpage{163}
\doi{10.3150/12-BEJ479} 

\makeatletter

\newcolumntype{d}[1]{D{.}{.}{#1}}

\newcommand{\mrmd}{\,\mathrm{d}}

\newcommand{\RMo}{\mathrm{o}}
\newcommand{\RMO}{\mathrm{O}}

\newtheorem{theorem}{Theorem}
\newtheorem{corollary}{Corollary}
\newtheorem{lemma}{Lemma}

\newremark{remark}{Remark}

\newcommand{\apgt}{\gtrsim}
\newcommand{\aplt}{\lesssim}

\newcommand{\vece}{{\mathbf{e}}}

\newcommand{\cT}{\mathcal{T}}

\newcommand{\bJ}{\mathbf{J}}
\newcommand{\rbA}{\mathbf{A}}
\newcommand{\rbB}{\mathbf{B}}
\newcommand{\rbH}{\mathbf{H}}
\newcommand{\rbI}{\mathbf{I}}
\newcommand{\rbV}{\mathbf{V}}
\newcommand{\rbR}{\mathbf{R}}
\newcommand{\sfT}{\bolds{\mathsf{T}}}
\newcommand{\sfU}{\bolds{\mathsf{U}}}
\newcommand{\sfX}{\bolds{\mathsf{X}}}
\newcommand{\sfZ}{\bolds{\mathsf{Z}}}

\newcommand{\bSigma}{\bolds{\Sigma}}

\newcommand{\bbeta}{\bolds{\beta}}
\newcommand{\bgamma}{\bolds{\gamma}}
\newcommand{\bvarphi}{\bolds{\varphi}}
\newcommand{\bDelta}{\bolds{\Delta}}
\newcommand{\vecpsi}{{\bolds{\varpsi}}}

\newcommand{\add}{{+}}

\newcommand{\bT}{\mathbf{T}}
\newcommand{\bB}{\mathbf{B}}
\newcommand{\bX}{\mathbf{X}}
\newcommand{\bY}{\mathbf{Y}}
\newcommand{\bZ}{\mathbf{Z}}
\newcommand{\bV}{\mathbf{V}}

\newcommand{\bbG}{\mathbb{G}}
\newcommand{\bbR}{\mathbb{R}}
\newcommand{\bbF}{\mathbb{F}}
\newcommand{\var}{\operatorname{var}}
\newcommand{\vecmu}{\bolds{\mu}}

\newcommand{\varpsi}{{{\psi}}}
\newcommand{\argmin}[0]{\arg\min}

\makeatother

\begin{document}
\begin{frontmatter}

\title{Efficient semiparametric estimation in generalized partially
linear additive models for longitudinal/clustered data}
\runtitle{Partially linear additive models for longitudinal
data}

\begin{aug}
\author[1]{\fnms{Guang} \snm{Cheng}\corref{}\thanksref{1}\ead[label=e1]{chengg@purdue.edu}},
\author[2]{\fnms{Lan} \snm{Zhou}\thanksref{2,e2}\ead[label=e2,mark]{lzhou@stat.tamu.edu}} \and
\author[2]{\fnms{Jianhua Z.} \snm{Huang}\thanksref{2,e3}\ead[label=e3,mark]{jianhua@stat.tamu.edu}}
\runauthor{G. Cheng, L. Zhou and J.Z. Huang} 
\address[1]{Purdue University, West Lafayette, IN 47907, USA.
\printead{e1}}
\address[2]{Texas A\&M University, College Station, TX 77843, USA.
\printead{e2};\\ \printead*{e3}}
\end{aug}

\received{\smonth{10} \syear{2011}}
\revised{\smonth{9} \syear{2012}}

%
\begin{abstract}
We consider efficient estimation of the Euclidean parameters in a
generalized partially linear additive models for
longitudinal/clustered data when multiple covariates need to be
modeled nonparametrically, and propose an estimation procedure
based on a spline approximation of the nonparametric part of the
model and the generalized estimating equations (GEE). Although the
model in consideration is natural and useful in many practical
applications, the literature on this model is very limited because
of challenges in dealing with dependent data for nonparametric
additive models. We show that the proposed estimators are
consistent and asymptotically normal even if the covariance
structure is misspecified. An explicit consistent estimate of the
asymptotic variance is also provided. Moreover, we derive the
semiparametric efficiency score and information bound under
general moment conditions. By showing that our estimators achieve
the semiparametric information bound, we effectively establish
their efficiency in a stronger sense than what is typically
considered for GEE. The derivation of our asymptotic results
relies heavily on the empirical processes tools that we develop
for the longitudinal/clustered data. Numerical results are used to
illustrate the finite sample performance of the proposed
estimators.
\end{abstract}

%

\end{frontmatter}

\section{Introduction}\label{sec1}

The partially linear model has become a widely used semiparametric regression
model because it provides a nice trade-off between model
interpretability and flexibility.
In a partially linear model, the mean of the outcome is assumed to
depend on some covariates $\bX$ parametrically and some other covariates
$\bT$ nonparametrically. Usually, the effects of~$\bX$ (e.g., treatment)
are of major interest, while the effects of $\bT$ (e.g., confounders)
are nuisance parameters.
Efficient estimation for partially linear models has been extensively
studied and well understood for independent data; see, for example,
Chen \cite{c88}, Speckman \cite{s88}, and Severini and Staniswalis
\cite{ss94}.
The book of H\"ardle, Liang and Gao \cite{hlg00} provides a
comprehensive review of the subject.

Efficient estimation of the Euclidean parameter (i.e., the parametric
component)
in the partially linear model for dependent data is by no means simple
due to
complication in data structure. Lin and Carroll \cite{lc01a,lc01b}
showed that,
whether a natural application of the local polynomial kernel
method can yield a semiparametric efficient estimator depends on
whether the covariate modeled nonparametrically is a cluster-level
covariate or not. Because the naive approach fails, Wang, Carroll and
Lin \cite{wcl05}
constructed a semiparametric efficient estimator by employing the
iterative kernel method of Wang \cite{w03} that can effectively
account for the within-cluster correlation. Alternatively,
Zhang \cite{z04}, Chen and Jin \cite{cj06}, and Huang, Zhang and Zhou
\cite{hzz07} constructed
semiparametric efficient estimators by extending the parametric
generalized estimating equations (GEE) of Liang and Zeger \cite{lz86}.
He, Zhu and Fung \cite{hzf02} and He, Fung and Zhu \cite{hfz05}
considered robust estimation, and
Leng, Zhang and Pan \cite{czp10} studied joint mean-covariance modeling
for the
partially linear model also by extending the GEE. In all these
development, only one covariate is modeled nonparametrically.

In many practical situations, it is desirable to model multiple
covariates nonparametrically. However, it is well known that
multivariate nonparametric estimation is subject to the curse of
dimensionality. A widely used approach for dimensionality
reduction is to consider an additive model for the nonparametric
part of the regression function in the partly linear model, which
in turn results in the partially linear additive model. Although
adapting this approach is a natural idea, there are major
challenges for estimating the additive model for dependent data.
Until only very recently, Carroll \textit{et al.} \cite{cmmy09} gave
the first
contribution on the partly linear additive model for
longitudinal/clustered data, focusing on a simple setup of the
problem, where there is the same number of observations per
subject/cluster, and the identity link function is used.

The goal of the paper is to give a thorough treatment of the
problem in the general setting that allows a monotonic link
function and unequal number of observations among
subjects/clusters. In this general setting, we derive the
semiparametric efficient score and efficiency bound to obtain a
benchmark for efficient estimation. In our derivation, we only
assume the conditional moment restrictions instead of any
distributional assumptions, for example, the multivariate Gaussian error
assumption employed in Carroll \textit{et al.} \cite{cmmy09}. It turns out
the definition
of the efficient score involves solving a system of complex
integral equations and there is no closed-form expression. This
fact rules out the feasibility of constructing efficient
estimators by plugging the estimated efficient influence function
into their asymptotic linear expansions.
We propose an estimation procedure that approximates the unknown
functions by splines and uses the generalized estimating
equations. To differentiate our procedure with the parametric
GEE, we refer to it as the extended GEE. We show that the extended
GEE estimators are semiparametric efficient if the covariance
structure is correctly specified and they are still consistent and
asymptotically normal even if the covariance structure is
misspecified. In addition, by taking advantage of the spline
approximation, we are able to give an explicit consistent estimate
of the asymptotic variance without solving the system of integral
equations that lead to the efficient scores. Having a closed-form
expression for the asymptotic variance is an attractive feature of
our method, in particular when there is no closed-form expression
of the semiparametric efficiency bound. Another attractive feature
of our method is the computational simplicity, there is no need to
resort to the computationally more demanding backfitting type
algorithm and numerical integration, as has been done in the
previous work on the same model.

As a side remark, one highlight of our mathematical rigor is
the careful derivation
of the smoothness conditions on the least favorable directions
from primitive conditions.
This rather technical but
important issue has not been well treated in the literature.
To develop the asymptotic theory in this paper, we rely heavily
on some new empirical process tools which we develop
by extending existing results from the i.i.d. case to the
longitudinal/clustered data.

The rest of the paper is organized as follows. Section \ref{sec2}
introduces the setup of the partially linear additive model and
the formulation of the extended GEE estimator. Section \ref{sec3} lists all
regularity conditions, derives the semiparametric efficient score and the
efficiency bound, and presents the asymptotic properties of the
extended GEE estimators.
Section \ref{sec4} illustrates the finite sample performance of the GEE
estimators using a simulation study and a real data. The proofs of some
nonasymptotic results and the sketched proofs of the main asymptotic
results are given in
the \hyperref[app]{Appendix}. The supplementary file discusses the properties of the least
favorable directions, presents the relevant empirical processes tools and
the complete proofs of all asymptotic results.

\textit{Notation.} For positive number sequences $a_n$ and $b_n$, let $a_n
\aplt b_n$ mean that $a_n/b_n$ is bounded, $a_n \asymp b_n$ mean that
$a_n \aplt b_n$ and $a_n \apgt b_n$, and $a_n\ll b_n$ mean that
$\lim_na_n/b_n=0$. For two positive semidefinite matrices $\rbA$ and
$\rbB$,
let $\rbA\geq\rbB$ mean that $\rbA-\rbB$ is positive semidefinite.
Define $x\vee y$ ($x\wedge y$) to be the maximum
(minimum) value of $x$ and $y$.
For any matrix $\rbV$, denote $\lambda_{V}^{\max}$ $(\lambda_{V}^{\min})$
as the largest (smallest) eigenvalue of $\rbV$.
Let $|\bV|$ denote the Euclidean norm of the vector $\bV$.
Let $\| a \|_{L_2}$ denote the usual
$L_2$ norm of a squared integrable function $a$, where the domain of
integration and the dominating measure should be clear from the context.

\section{The model setup}\label{sec2}
Suppose that the data consist of $n$ clusters with the $i$th
($i=1,\ldots, n$) cluster having $m_i$ observations. In
particular, for longitudinal data a cluster represents an
individual subject. The data from different clusters are
independent, but correlation may exist within a cluster. Let
$Y_{ij}$ and $(\bX_{ij}, \bT_{ij})$ be the response variable and
covariates for the $j$th ($j=1,\ldots, m_i$) observation in the
$i$th cluster. Here $\bX_{ij} = (X_{ij1},\ldots, X_{ijK})'$ is a
$K \times1$ vector and $\bT_{ij} = (T_{ij1},\ldots, T_{ijD})'$ is
a $D \times1$ vector.
We consider the marginal model
%
%
\begin{equation}
\label{eqmodel1} \mu_{ij} = E(Y_{ij} | \bX_{ij},
\bT_{ij}),
\end{equation}
and the marginal mean $\mu_{ij}$ depends on covariates $\bX_{ij}$
and $\bT_{ij}$ through a known monotonic and differentiable link
function $\mu(\cdot)$:
%
\begin{eqnarray}
\label{eqmodel2} \mu_{ij} & = & \mu\bigl(\bX_{ij}'
\bbeta+ \theta_\add(\bT_{ij})\bigr)
\nonumber\\[-8pt]\\[-8pt]
& = & \mu\bigl(\bX_{ij}'\bbeta+ \theta_1(T_{ij1})
+ \cdots+ \theta_D(T_{ijD})\bigr),\nonumber
\end{eqnarray}
where $\bbeta$ is a $K \times1$ vector, and $\theta_\add(\mathbf
{t})$ is
an additive function with $D$ smooth additive component functions
$\theta_d(t_d)$, $1\leq d\leq D$. For\vspace*{1pt} identifiability, it is
assumed that $\int_{\cT_d} \theta_d(t_d) \mrmd t_d =0$, where $\cT_d$
is the compact support of the covariate $T_{ijd}$. Applications of
marginal models for longitudinal/clustered data are common in the
literature (Diggle \textit{et al.} \cite{dhlz02}).\looseness=-1

Denote
\begin{eqnarray*}
\mathbf{Y}_i &=& %
\pmatrix{ Y_{i1}
\cr
\vdots
\cr
Y_{im_i} },\qquad \vecmu_i= %
\pmatrix{
\mu_{i1}
\cr
\vdots
\cr
\mu_{im_i} },\qquad \sfX_i =
\pmatrix{ \bX_{i1}'
\cr
\vdots
\cr
\bX_{im_i}' },\qquad \sfT_i = %
\pmatrix{
\bT_{i1}'
\cr
\vdots
\cr
\bT_{im_i}'
},
\\
\theta_\add(\sfT_i) &=& %
\pmatrix{
\theta_\add(\bT_{i1})
\cr
\vdots
\cr
\theta_\add(
\bT_{im_i}) },\qquad \mu\bigl(\sfX_i\bbeta+ \theta_\add(
\sfT_i)\bigr)= %
\pmatrix{ \mu\bigl(\bX_{i1}'
\bbeta+\theta_\add(\bT_{i1})\bigr)
\cr
\vdots
\cr
\mu\bigl(
\bX_{im_i}'\bbeta+\theta_\add(\bT_{im_i})
\bigr) }.
\end{eqnarray*}
Here and hereafter, we make the notational convention that
application of a multivariate function to a matrix is understood
as application to each row of the matrix, and similarly
application of a univariate function to a vector is understood as
application to each element of the vector. Using matrix notation,
our model representation (\ref{eqmodel1}) and (\ref{eqmodel2})
can be written as
%
\begin{equation}
\label{eqmodel3} \vecmu_i= E(\bY_i|
\sfX_i, \sfT_i) = \mu\bigl(\sfX_i\bbeta+
\theta_\add(\sfT_i)\bigr).
\end{equation}


Note that in our modeling framework no distributional assumptions
are imposed on the data other than the moment conditions specified
in (\ref{eqmodel1}) and (\ref{eqmodel2}). In particular,
$\sfX_i$ and $\sfT_i$ are allowed to be dependent, as commonly
seen for longitudinal/clustered data. Let $\bSigma_i =
\operatorname{var}(\bY_i|\sfX_i, \sfT_i)$ be the true covariance matrix
of $\bY_i$. Following the generalized estimating equations (GEE)
approach of Liang and Zeger \cite{lz86}, we introduce a working covariance
matrix $\rbV_i = \rbV_i(\sfX_i, \sfT_i)$ of $\bY_i$, which can
depend on a nuisance finite-dimensional parameter vector $\tau$
distinct from $\bbeta$. In the parametric setting, Liang and Zeger
\cite{lz86}
showed that, consistency of the GEE estimator is guaranteed even
when the covariance matrices are misspecified, and estimation
efficiency will be achieved when the working covariance matrices
coincide with the true covariance matrices, that is, when $
\rbV_i(\tau^*)=\bSigma_i$ for some $\tau^*$. In this paper, we
shall establish a similar result in a semiparametric context.

To estimate the functional parameters, we use basis approximations
(e.g., Huang, Wu and Zhou~\cite{hwz02}). We approximate each component function
$\theta_d(t_d)$ of the additive function $\theta_{\add}(\mathbf
{t})$ in~(\ref{eqmodel2}) by a basis expansion, that is,
%
\begin{equation}
\label{eqspl-app} \theta_d(t_d) \approx\sum
_{q=1}^{Q_d} \gamma_{dq} B_{dq}(t_d)=
\bB'_d(t_d) \bgamma_d,
\end{equation}
where $B_{dq}(\cdot), q=1,\ldots, Q_d$, is a system of basis
functions, which is denoted as a vector $\bB_d(\cdot)=
(B_{d1}(\cdot),\ldots, B_{dQ_d}(\cdot))'$, and $\bgamma_d=
(\gamma_{d1},\ldots, \gamma_{dQ_d})'$ is a vector of
coefficients.\vadjust{\goodbreak}
In principle, any basis system can be used, but B-splines are used
in this paper for their good approximation properties. In fact, if
$\theta_d(\cdot)$ is continuous, the spline approximation can be
chosen to satisfy $\sup_t|\theta_d(t)- \bB_d'(t)\bgamma_d| \to0$
as $Q_d \to\infty$, and the rate of convergence can be
characterized based on the smoothness of $\theta_d(\cdot)$; see
de Boor \cite{d01}.

It follows from (\ref{eqspl-app}) that
%
\begin{equation}
\label{eqspl-app2} \theta_\add(\bT_{ij}) \approx\sum
_{d=1}^D \sum_{q=1}^{Q_d}
\gamma_{dq} B_{dq}(T_{ijd}) = \sum
_{d=1}^D \bB'_d(T_{ijd})
\bgamma_d = \bZ_{ij}'\bgamma,
\end{equation}
where $\bZ_{ij} = (\bB'_1(T_{ij1}),\ldots, \bB'_D(T_{ijD}))'$, and
$\bgamma= ({\bgamma_1}',\ldots, {\bgamma_D}')'$. Denoting $\sfZ_i
= (\bZ_{i1},\ldots,\allowbreak \bZ_{im_i})'$, (\ref{eqmodel3}) and
(\ref{eqspl-app2}) together imply that
%
\begin{equation}
\label{eqmodel4} \vecmu_i= E(\bY_i|
\sfX_i, \sfT_i) \approx\mu(\sfX_i\bbeta+
\sfZ_i \bgamma).
\end{equation}
Thus, the Euclidean parameters and functional parameters are
estimated jointly by minimizing the following weighted least
squares criterion
%
\begin{equation}
\label{eqwls} \sum_{i=1}^n \bigl\{
\bY_i - \mu(\sfX_i\bbeta+ \sfZ_i\bgamma)
\bigr\}' \rbV_i^{-1} \bigl\{
\bY_i - \mu(\sfX_i\bbeta+ \sfZ_i\bgamma)
\bigr\}
\end{equation}
or, equivalently, by solving the estimating equations
%
\begin{equation}
\label{eqgee1} \sum_{i=1}^n
\sfX_i' \bDelta_i \rbV_i^{-1}
\bigl\{\bY_i - \mu(\sfX_i\bbeta+ \sfZ_i
\bgamma)\bigr\} = 0
\end{equation}
and
%
\begin{equation}
\label{eqgee2} \sum_{i=1}^n
\sfZ_i' \bDelta_i \rbV_i^{-1}
\bigl\{\bY_i - \mu(\sfX_i\bbeta+ \sfZ_i
\bgamma)\bigr\} = 0,
\end{equation}
where\vspace*{1pt} $\bDelta_i$ is a diagonal matrix with the diagonal elements
being the first derivative of $\mu(\cdot)$ evaluated at
$\bX_{ij}'\bbeta+ \bZ_{ij}'\bgamma$, $j=1,\ldots, m_i$. Denoting
the minimizer of (\ref{eqwls}) as $\widehat\bbeta$ and
$\widehat\bgamma$, then $\widehat\bbeta$ estimates the parametric
part of the model, and $\widehat\theta_1(\cdot) =
\bB'_1(\cdot)\widehat\bgamma_1,\ldots,\widehat\theta_D(\cdot
) =
\bB'_D(\cdot)\widehat\bgamma_D$ estimate the nonparametric part of
the model. We refer to these estimators the extended GEE
estimators. In this paper, we shall show that, under regularity
conditions, $\widehat\bbeta$ is asymptotically normal and, if the
correct covariance structure is specified, it is semiparametric
efficient, and also show that $\widehat\theta_{d}(\cdot)$ is a
consistent estimator of the true nonparametric function
$\theta_{d}(\cdot)$, $d=1,\ldots, D$.

When the link function $\mu(\cdot)$ is the identity function, the
minimizer of the weighted least squares (\ref{eqwls}) or the
solution to the estimating equations (\ref{eqgee1}) and
(\ref{eqgee2}) has a closed-from expression:
\[
\pmatrix{ \widehat\bbeta
\cr
\widehat\bgamma} %
= \Biggl(
\sum_{i=1}^{n}\sfU_i'
\rbV_i^{-1}\sfU_i \Biggr)^{-1}\sum
_{i=1}^{n}\sfU_i'
\rbV_i^{-1}\bY_i,
\]
where $\sfU_i=(\sfX_i, \sfZ_i)$.\vadjust{\goodbreak}

\section{Theoretical studies of extended GEE estimators}\label{sec3}
\subsection{Regularity conditions}\label{sec3.1}\label{subassump}
We state the regularity conditions needed for the
theoretical results in this paper. For the asymptotic
analysis, we assume that the number of
individuals/clusters goes to infinity while the number of
observations per individual/cluster remains bounded.

\begin{enumerate}[C1.]
\item[C1.] The random variables $T_{ijd}$ are bounded, uniformly
in $i=1,\ldots, n$, $j=1,\ldots, m_i$ and $d=1,\ldots, D$. The
joint distribution of any pair of $T_{ijd}$ and $T_{ij'd'}$ has a
density $f_{ijj'dd'}(t_{ijd},t_{ij'd'})$ with respect to the
Lebesgue measure. We assume that $f_{ijj'dd'}(\cdot, \cdot)$ is
bounded away from 0 and infinity, uniformly in $i=1,\ldots, n$,
$j, j'=1,\ldots, m_i$, and $d,d'=1,\ldots,D$.

\item[C2.] The first covariate is constant 1, that is, $X_{ij1}
\equiv1$. The random variables $X_{ijk}$ are bounded, uniformly
in $i=1,\ldots, n$, $j=1,\ldots, m_i$ and $k=2,\ldots, K$. The
eigenvalues of $E\{\bX_{ij} \bX_{ij}'|\bT_{ij}\}$ are bounded away
from 0, uniformly in $i=1,\ldots, n$, $j=1,\ldots, m_i$.

\item[C3.] The eigenvalues of true covariance matrices $\bSigma_i$
are bounded away from 0 and infinity, uniformly in $i=1,\ldots,
n$.

\item[C4.] The eigenvalues of the working covariance matrices
$\rbV_i$ are bounded away from 0 and infinity, uniformly in $i=1,\ldots, n$.
\end{enumerate}

Conditions similar to C1--C4 were used and discussed in
Huang, Zhang and Zhou \cite{hzz07} when considering partially linear
models with the
identity link. Condition C1 is also used to ensure identifiability
of the additive components, see Lemma 3.1 of Stone \cite{s94}.
Condition C1 implies that the marginal density $f_{ijd}(\cdot)$ of
$T_{ijd}$ is bounded away from 0 on its support, uniformly in
$i=1,\ldots, n$, $j=1,\ldots, m_i$, and $d=1,\ldots,D$. The
condition on eigenvalues in C2 prevents the multicollinearity of
the covariate vector $\bX_{ij}$ and ensures the identifiability of
$\bbeta$. Since we assume that the cluster size (or the number of
observations per subject) is bounded, we expect C3 is in general
satisfied. Note that a zero eigenvalue of $\Sigma_i$ indicates that
there is a perfect linear relation among the residuals from subject
$i$, which is unlikely to happen in reality.

Denote the true values of $\bbeta$ and $\theta_\add(t)$ by
$\bbeta_0$ and $\theta_{0,\add}(t)$, respectively.

\begin{enumerate}[C5.]
\item[C5.] (i) The link function $\mu$ is strictly monotone and
has continuous second derivative; (ii)~$\inf_s \mu'(s) = c_1
> 0$; (iii) $\mu'$ and $\mu''$ are locally bounded around
$\mathbf{x}^T\bbeta_0 +\theta_{0,\add}(\mathbf{t})$; (iv)~$\mu
(\pm v)$ increases
slower than $v^L$ as $v \to\infty$ for some $L>0$.
\end{enumerate}

Denote $e_{ij} = Y_{ij} - \mu_{ij}$ and $\vece_i = (e_{i1},\ldots,
e_{im_i})'$.

\begin{enumerate}[C6.]
\item[C6.] The errors are uniformly sub-Gaussian, that is,
%
\begin{equation}
\label{eqgl-assumption2} \max_{i=1,\ldots,n} M_0^2 E
\bigl\{\exp\bigl(|\vece_i|^2/M_0^2
\bigr)-1 \vert\sfX_i, \sfT_i\bigr\} \leq
\sigma_0^2\qquad \forall n\mbox{, a.s.}
\end{equation}
for some fixed positive constants $M_0$ and $\sigma_0$.
\end{enumerate}

Condition C5 on the link function is satisfied in all practical
situations. The sub-Gaussian condition C6 relaxes the strict
multivariate Gaussian error assumption, and is commonly used in
the literature when applying the empirical process theory.

For $i=1,\ldots, n$, let $\bDelta_{i0}$ be a diagonal matrix with
the $j$th diagonal element being the first derivative of
$\mu(\cdot)$ evaluated at $\bX_{ij}'\bbeta_0 +
\theta_{0,\add}(\bT_{ij})$, $j=1,\ldots, m_i$. Let $\bX_{ik}$
denote the $k$th column of the matrix $\sfX_i$. For any additive
function $\varphi_\add(\mathbf{t}) = \varphi_1(t_1) + \cdots+
\varphi_D(t_D)$, $\mathbf{t}= (t_1,\ldots, t_D)'$, denote
$\varphi_\add(\sfT_i)=(\varphi_\add(\bT_{i1}),\ldots,\varphi_\add(\bT
_{im_i}))'$.
Let $\varphi^*_{k,\add}(\cdot)$ be the additive function
$\varphi_{k,\add}(\cdot)$ that minimizes
%
\begin{equation}\label{varphi}
\sum_{i=1}^n E \bigl[\bigl\{
\bX_{ik}- \varphi_{k,\add}(\sfT_i)\bigr
\}'\bDelta_{i0} \rbV_i^{-1}
\bDelta_{i0}\bigl\{\bX_{ik}- \varphi_{k,\add}(
\sfT_i)\bigr\} \bigr].
\end{equation}
Denote $\bvarphi^\ast_\add(\sfT_i)=(\varphi_{1,\add}^\ast(\sfT_i),\ldots,\varphi_{K,\add}^\ast(\sfT_{i}))$ and define
\[
\rbI_V \equiv\lim_{n}\frac{1}{n}\sum
_{i=1}^n E \bigl[ \bigl\{\sfX_{i}-
\bvarphi_\add^\ast(\sfT_i)\bigr\}'
\bDelta_{i0} \rbV_i^{-1}\bDelta_{i0}
\bigl\{\sfX_{i}-\bvarphi^\ast_\add(
\sfT_i)\bigr\} \bigr].
\]

\begin{enumerate}[C7.]
\item[C7.] The matrix $\rbI_V$ is positive definite.
\end{enumerate}

Condition C7 is a positive information requirement that ensures
the Euclidean parameter $\bbeta$ can be root-$n$ consistently
estimated. When $\rbV_i$ is specified to be the true covariance
matrix $\bSigma_i$ for all $i$, $\varphi_{k,\add}^\ast(\cdot)$
reduces to the least favorable direction
$\psi_{k,\add}^\ast(\cdot)$ in the definition of efficient score
function and $\rbI_V$ reduces to the efficient information matrix
$\rbI_{\mathrm{eff}}$; see Section \ref{secscore}.

For $d =1,\ldots, D$, let $\bbG_d= \{\bB_d'(t)\bgamma_d\}$ be a
linear space of splines with degree $r$ defined on the support
$\cT_d$ of $T_{ijd}$. Let $\bbG_\add= \bbG_1 + \cdots+ \bbG_D$ be
the additive spline space. We allow the dimension of $\bbG_d$,
$1\leq d\leq D$, and $\bbG_\add$ to depend on $n$, but such
dependence is suppressed in our notation to avoid clutter. For
each spline space, we require that the knot sequence satisfies the
quasi-uniform condition, that is, $\max_{j, j'} (u_{n,j+r+1} -
u_{n,j})/(u_{n,j'+r+1} -u_{n,j'})$ is bounded uniformly in $n$ for
knots $\{u_{n, j}\}$. Let
\[
\rho_n = \max\Bigl\{ \inf_{g\in\bbG_\add}\bigl\|g(\cdot)-
\theta_{0,\add}(\cdot)\bigr\|_\infty, \max_{1\leq k\leq K}
\inf_{g\in\bbG_\add} \bigl\|g(\cdot)- \varphi_{k,\add}^*(\cdot)\bigr\|_\infty
\Bigr\}
\]
and $Q_n = \max\{Q_d= \operatorname{dim}(\bbG_d), 1\le d \le D\}$.

\begin{enumerate}[C8.]
\item[C8.] (i) $\lim_n Q_n^2 \log^4 n /n =0$, (ii) $\lim_n n
\rho_n^4=0$.
\end{enumerate}

Condition C8(i) characterizes the growth rate of the dimension of
the spline spaces relative to the sample size. Condition C8(ii)
describes the requirement on the best rate of convergence that the
functions $\theta_{0,\add}(\cdot)$ and
$\varphi^*_{k,\add}(\cdot)$'s can be approximated by functions in
the spline spaces. These requirements can be quantified by
smoothness conditions on $\theta_{0,\add}(\cdot)$ and
$\varphi^*_{k,\add}(\cdot)$'s, as follows. For $\alpha>0$, write
$\alpha= \alpha_0 + \alpha_1$, where $\alpha_0$ is an integer and
$0<\alpha_1 \leq1$. We say a function is $\alpha$-smooth, if its
derivative of order $\alpha_0$ satisfies a H\"older condition with
exponent $\alpha_1$. If all additive components of
$\theta_{0,\add}(\cdot)$ and $\varphi_{k,\add}^*(\cdot)$'s are
$\alpha$-smooth, and the degree $r$ of the splines satisfies $r
\geq\alpha-1$, then, by a standard result from approximation
theory, $\rho_n \asymp Q_n^{-\alpha}$ for $\alpha>1/2$
(Schumaker \cite{s81}). Condition C8 thus can be replaced by the following
condition.

\begin{enumerate}[C$8'.$]
\item[C$8'.$] (i) $\lim_n Q_n^2 \log^4 n /n =0$; (ii) additive
components of $\theta_{0,\add}(\cdot)$ and
$\varphi_{k,\add}^\ast(\cdot)$, $k=1,2,\ldots,\allowbreak K$, are
$\alpha$-smooth for some $\alpha>1/2$; (iii) $\lim_n
Q_n^{4\alpha}/n=\infty$.
\end{enumerate}

Since $\varphi_{k,\add}^*$ is only implicitly defined, it is
important to verify its smoothness requirement from primitive
conditions. In the supplementary file (Cheng, Zhou and Huang \cite
{chz12}), that is, Section S.1, we shall show that
$\varphi_{k,\add}^*(\cdot)$ solves a system of integral equations
and its smoothness is implied by smoothness requirements on the
joint density of $\sfX_i$ and $\sfT_i$.

\subsection{Semiparametric efficient score and efficiency
bound}\label{sec3.2} \label{secscore}
For estimating the Euclidean parameter in a semiparametric model,
the efficiency bound provides a useful benchmark for the optimal
asymptotic behaviors (e.g., Bickel \textit{et al.} \cite{bkrw93}). In this
subsection, we
give the semiparametric efficient score and efficient information
matrix when the covariance structure is correctly specified. We do not
make the
normality assumption on the error distribution in the derivations.

The models studied in this paper have more than one nuisance function
so that the efficient score function for $\bbeta$, denoted as
$\ell^\ast_{\beta}$, is obtained by projecting onto a sum-space. In
Lem\-ma~\ref{lemmaeff-score} below,
we construct $\ell^\ast_{\beta}$ by the two-stage projection approach
(Sasieni \cite{s92}). Recall that $\vece_i={\bY}_i-\mu({\sfX
}_i\bbeta
+\theta_\add({\sfT}_i))$,
where
$\theta_\add(\mathbf{t}) = \theta_1(t_1) + \cdots+ \theta_D(t_D)$.
Write $f_i({\bolds{\mathsf{x}}}_i, {\bolds{\mathsf{t}}}_i,
{{\mathbf{y}}}_i-\mu({\bolds{\mathsf{x}}}_i\bbeta+\theta_\add
({\bolds{\mathsf{t}}}_i)))$
as the joint density of $({\sfX}_i, {\sfT}_i, {\bY}_i)$ for the $i$th
cluster. We assume that $f_i(\cdot,\cdot,\cdot)$ is smooth, bounded
and satisfies $\lim_{|e_{ij}|\rightarrow\infty}f_i(\cdot,\cdot
,{\mathbf
e}_i)=0$
for all $j=1,\ldots,m_i$.
%
\begin{lemma}\label{lemmaeff-score}
The efficient score has the expression
$\ell^\ast_\beta=
(\ell^\ast_{\beta,1},\ldots, \ell^\ast_{\beta,K})'$ with
%
\begin{equation}
\label{eqeffsco} \ell^\ast_{\beta, k} = \sum
_{i=1}^n \bigl(\bX_{ik}-
\psi^*_{k,\add
}(\sfT_i) \bigr)'
\bDelta_{i0}\bSigma_i^{-1}\bigl[{
\bY}_i -\mu\bigl({\sfX}_i\bbeta_0+
\theta_{0,\add}({\sfT}_i)\bigr)\bigr],
\end{equation}
where $\psi^*_{k,\add}(\mathbf{t}) = \sum_{d=1}^{D}\psi_{kd}^\ast
(t_d)$ satisfies
%
\begin{equation}
\label{psikast} \sum_{i=1}^n E \bigl[
\bigl({\bX}_{ik}- \psi_{k,\add}^\ast(\sfT)
\bigr)' \bDelta_{i0}\bSigma_i^{-1}
\bDelta_{i0}\psi_d({\bT}_{id}) \bigr] =0
\end{equation}
for any $\psi_{d}(t_d)\in L_2(\mathcal{T}_d)$, $d=1,\ldots, D$.
\end{lemma}

The form of $\ell^\ast_{\beta,k}$ when $D=1$ coincides with that
derived in the partially linear models, for example, Lin and Carroll \cite
{lc01a}, under the strict multivariate Gaussian error assumption.
In the supplementary material (Cheng, Zhou and Huang \cite{chz12}), we
shall see that $\psi_{kd}^\ast(t_d)$'s (or, more generally, $\varphi
_{kd}^\ast(t_d)$'s) solve a Fredholm integral equation of the second kind
(Kress \cite{k99}), and do not have a closed-form expression. In the
same file, we also show that $\psi_{kd}^\ast(t_d)$'s (or, more
generally, $\varphi_{kd}^\ast(t_d)$'s)
are well defined and have nice properties such as boundedness and
smoothness under reasonable assumptions on the joint density of
$\sfX_i$ and $\sfT_i$. These properties are crucial for the
feasibility to construct semiparametric efficient estimators but
are not carefully studied in the literature.

The semiparametric efficient information matrix for $\bbeta$ is
%
\begin{eqnarray}
\label{eqinfo} \rbI_{\mathrm{eff}} & \equiv &\lim_n \frac{1}{n} E
\bigl(\ell_\beta^* \ell_\beta^{*\prime}\bigr)
\nonumber\\[-8pt]\\[-8pt]
& = &\lim_n \frac{1}{n} \sum_{i=1}^n
E \bigl[\bigl\{\sfX_i- \vecpsi^*_{\add}(\sfT_i)
\bigr\}' \bDelta_{i0} \bSigma_i^{-1}
\bDelta_{i0} \bigl\{\sfX_i- \vecpsi^*_{\add}(
\sfT_i)\bigr\} \bigr],\nonumber
\end{eqnarray}
where $\vecpsi_{\add}^*(\sfT_i) = (\psi_{1,\add}^*(\sfT_i),\ldots, \psi_{K,\add}^*(\sfT_i))$. The efficient information matrix
$\rbI_{\mathrm{eff}}$ here is the
same as the quantity $\rbI_V$ in condition C7 when $\rbV_i = \bSigma_i$.
In the above result, different subjects/clusters need not have
the same number of observations and thus $(\sfX_i, \sfT_i)$
may not be identically distributed.
In the special case that $(\sfX_i, \sfT_i)$ are i.i.d., the efficient
information can be simplified to
\[
\rbI_{\mathrm{eff}}= E\bigl[\bigl\{\sfX_i- \vecpsi^*_\add(
\sfT_i)\bigr\}' \bDelta_{i0}
\bSigma_i^{-1} \bDelta_{i0} \bigl\{
\sfX_i- \vecpsi^*_\add(\sfT_i)\bigr\}\bigr],
\]
where the $k$th component of $\vecpsi^*_\add$ satisfies
\[
E \bigl[\bigl\{{\bX}_{ik}- \psi_{k,\add}^\ast(\sfT)
\bigr\}' \bDelta_{i0}\bSigma_i^{-1}
\bDelta_{i0}\psi_d({\bT}_{id}) \bigr] =0
\]
for any $\psi_{d}(t_d)\in L_2(\mathcal{T}_d)$, $d=1,\ldots, D$.

The function $\vecpsi_{\add}^*(\sfT_i)$ involved in the efficient
information matrix (\ref{eqinfo}) actually corresponds to the
least favorable direction (LFD) along $\theta_{0,\add}(\sfT_i)$ in
the least favorable submodel (LFS). To provide an intuitive
interpretation, we assume for simplicity that
$f_i(\vece_i|\bolds{\mathsf{x}}_i,\bolds{\mathsf
{t}}_i)\sim N(0,\bSigma_i)$. Given the above
distributional assumption, the parametric submodel (indexed by
$\varepsilon$) passing through $(\bbeta_0, \theta_{0,\add})$ is
constructed as
%
\begin{equation}\label{lfs}
\varepsilon\mapsto-\frac{1}{2}\sum_{i=1}^{n}
\bigl[ {{\mathbf{y}}}_i-\vecmu_i(\varepsilon)
\bigr]'\bSigma_i^{-1}\bigl[ {{
\mathbf{y}}}_i-\vecmu_i(\varepsilon)\bigr],
\end{equation}
where $\vecmu_i(\varepsilon)=\mu\{{\bolds{\mathsf
{x}}}_i(\bbeta_0+\varepsilon{\mathbf{v}}
)+[\theta_{0,\add}({\bolds{\mathsf{t}}}_i)+\varepsilon h_{\add
}(\bolds{\mathsf{t}}_i)]\}$, for some
vector ${\mathbf{v}}\in\mathbb{R}^{K}$ and perturbation
direction $h_{\add}(\cdot
)$ around $\theta_{0,\add}(\cdot)$. For any fixed ${\mathbf
{v}}$, the
information matrix for the parametric submodel (evaluated at
$\varepsilon=0$) is calculated as
\[
\rbI_{\mathrm{para}}(h_{\add})=\lim_n \frac{1}{n} \sum
_{i=1}^n E \bigl[\bigl\{
\sfX_i{\mathbf{v}}+ h_{\add}(\sfT_i)\bigr
\}' \bDelta_{i0} \bSigma_i^{-1}
\bDelta_{i0} \bigl\{\sfX_i{\mathbf{v}}+
h_{\add}(\sfT_i)\bigr\} \bigr].
\]
The minimum $\rbI_{\mathrm{para}}(h_{\add})$ over all possible
perturbation directions is known as the semiparametric efficient
information for ${\mathbf{v}}'\bbeta$ (Bickel \textit{et al.} \cite
{bkrw93}). The
parametric submodel
achieving the minimum is called the LFS and the associated direction is
called LFD. By calculating the Fr\'{e}chet derivative of the quadratic
function $h_\add\mapsto\rbI_{\mathrm{para}}(h_{\add})$ and considering
(\ref{psikast}), we can easily show that its minimum is achieved when
$h_\add=-\vecpsi^*_\add{\mathbf{v}}$. In view of the above
discussion, the
efficient information for $\bbeta$ becomes the $\rbI_{\mathrm{eff}}$ defined
in (\ref{eqinfo}).
%
\begin{remark}
Our derivation of the efficient score and efficient information
matrix also applies when $\bT$ is a cluster level covariate, that is,
$T_{ijd}=T_{id}$ for $j=1,\ldots,m_i$, $d=1,\ldots, D$. Let
$\widetilde\bT_i=(T_{i1},\ldots,T_{iD})'$. In this case, we only
need to replace $\varpsi_{k,\add}^\ast(\sfT_i)$ and
$\varpsi_{\add}^\ast(\sfT_i)$ by $\varpsi_{k,\add}^\ast
(\widetilde
\bT_i)\mathbf{1}$ and $\mathbf{1}(\varpsi^\ast_{1,\add}(\widetilde
\bT_i),\ldots,\varpsi^\ast_{K,\add}(\widetilde\bT_i))$, where
$\mathbf{1}$ is an $m_i$-vector of ones, and do similar changes for
$\varpsi_{k,\add}(\sfT_i)$ and $\varpsi_{\add}(\sfT_i)$. It is
interesting to note that, when $(\bY_i,\sfX_i, \sfT_i)$ are
i.i.d., then $\varpsi_{k,\add}^\ast(\cdot)$ has a closed form
expression:
\[
\varpsi_{k,\add}^\ast(\mathbf{t})=\frac{E(\bX_{ik}'\bDelta_{i0}\bSigma_i^{-1}
\bDelta_{i0}\mathbf{1}|\widetilde\bT_i=\mathbf{t})}{E(
\mathbf{1}'\bDelta_{i0}\bSigma_i^{-1}
\bDelta_{i0}\mathbf{1}|\widetilde\bT_i=\mathbf{t})}.
\]
\end{remark}

\subsection{Asymptotic properties}\label{sec3.3}\label{secasypro}
In this subsection, we assume that the dimension of the Euclidean
parameter, that is, $K$, is fixed.
Define $f_0(\mathbf{x},\mathbf{t}) = \mu(\mathbf{x}'\bbeta_0 +
\theta_{0,\add}(\mathbf{t}))$.
Define
\[
\bbF_n = \bigl\{f(\mathbf{x}, \mathbf{t})\dvtx  f(\mathbf{x},\mathbf{t})
= \mu\bigl(\mathbf{x}'\bbeta+ g(\mathbf{t})\bigr), \beta\in
\bbR^K, g\in\bbG_\add\bigr\}.
\]
The extended GEE estimator can be written as
\[
\arg\min_{f \in\bbF_n} \frac{1}{n} \sum_{i=1}^n
\bigl\{\bY_i - f(\sfX_i, \sfT_i)\bigr
\}' \rbV_i^{-1} \bigl\{\bY_i -
f(\sfX_i, \sfT_i)\bigr\}.
\]
The minimizer is $\widehat f_n(\mathbf{x},\mathbf{t}) = \mu(\mathbf
{x}^T\widehat{\bbeta}_V +
\widehat\theta(\mathbf{t}))$ where $\widehat\theta(\mathbf{t}) =
\bB'(\mathbf{t})\widehat
\bgamma$.
The subscript of $\widehat\bbeta_V$ denotes
the dependence on the working covariance matrices.

According to condition C8 (or C$8'$), there is an additive
spline function $\theta_n^*(\mathbf{t}) = \bB'(\mathbf{t})\bgamma^*\in
\bbG_\add$
such that $\|\theta_n^* - \theta_{0,\add} \|_\infty\lesssim\rho_n
\to0$.
Then $f_n^*(\mathbf{x}, \mathbf{t}) = \mu(\mathbf{x}'\bbeta_0 +
\theta_n^*(\mathbf{t}))$ is
a spline-based approximation to the regression function.
Define
\[
\langle\xi_1, \xi_2\rangle_n =
\frac{1}{n} \sum_i \xi_1'(
\sfX_i, \sfT_i) \rbV_i^{-1}
\xi_2(\sfX_i, \sfT_i)
\]
and $\|\xi\|_n^2 = \langle\xi, \xi\rangle_n$.
%
\begin{theorem}[(Consistency)]\label{thmrate} The following results hold:
%
\begin{eqnarray}
\label{convrate}
\bigl\|\widehat f_n - f_n^*\bigr\|_n^2
&=& \RMO_P \bigl(Q_n\log^2 n/n\vee
\rho_n^2 \bigr),
\\
\label{noncons}
\bigl\|\widehat f_n - f_n^*\bigr\|_\infty
& = &
\RMo_P(1),
\\
\label{eqconsis-whole}
\|\widehat{f}_n - f_0 \|_\infty
& = &
\RMo_P(1),
\\
\label{eqconsis}
\widehat\bbeta_V \stackrel{P} {\rightarrow} \bbeta_0,\qquad
\|\widehat\theta- \theta_{0,\add}\|_\infty&=&
\RMo_P(1).
\end{eqnarray}
\end{theorem}

Theorem \ref{thmrate} says that the extended GEE estimators are consistent
in estimating the parametric and nonparametric components of the model.
Next we show that, our extended GEE estimator $\widehat{\bbeta}$ is
asymptotically normal even when the working covariance matrices
$\rbV_i$'s are not necessarily the same as the true ones.

Denote $\sfU_i = (\sfX_i, \sfZ_i)$ as before. Let
%
\begin{eqnarray}
\label{eqH} \rbH& = &\sum_{i=1}^n
\sfU_i' \bDelta_{i0} \rbV_i^{-1}
\bDelta_{i0} \sfU_i \equiv%
\pmatrix{
\rbH_{11} & \rbH_{12}
\cr
\rbH_{21} &
\rbH_{22} }
\nonumber\\[-8pt]\\[-8pt]
& = &
\pmatrix{ \displaystyle \sum_{i=1}^n
\sfX_i' \bDelta_{i0} \rbV_i^{-1}
\bDelta_{i0} \sfX_i & \displaystyle \sum_{i=1}^n
\sfX_i' \bDelta_{i0} \rbV_i^{-1}
\bDelta_{i0} \sfZ_i
\cr
\displaystyle \sum_{i=1}^n
\sfZ_i' \bDelta_{i0} \rbV_i^{-1}
\bDelta_{i0} \sfX_i & \displaystyle \sum_{i=1}^n
\sfZ_i' \bDelta_{i0} \rbV_i^{-1}
\bDelta_{i0} \sfZ_i }.\nonumber
\end{eqnarray}
By the block matrix form of matrix inverse,
%
\begin{eqnarray}
\label{eqHinv} \pmatrix{ \rbH_{11} & \rbH_{12}
\cr
\rbH_{21} & \rbH_{22} } %
^{-1} & = &
\pmatrix{ \rbH^{11} & \rbH^{12}
\cr
\rbH^{21} & \rbH^{22} } %
\nonumber\\[-8pt]\\[-8pt]
& = &
\pmatrix{ \rbH^{-1}_{11\cdot2} & -
\rbH_{11 \cdot2}^{-1} \rbH_{12} \rbH_{22}^{-1}
\vspace*{1pt}\cr
- \rbH_{22 \cdot1}^{-1} \rbH_{21} \rbH_{11}^{-1}
& \rbH_{22\cdot1}^{-1} },\nonumber
\end{eqnarray}
where $\rbH_{11\cdot2}= \rbH_{11} - \rbH_{12}\rbH^{-1}_{22} \rbH_{21}$ and
$\rbH_{22\cdot1}= \rbH_{22} - \rbH_{21}\rbH^{-1}_{11} \rbH_{12}$.
Define
\[
\rbR^{\vartriangle}(\widehat\bbeta_V) \equiv\rbH^{11}
\sum_{i=1}^n \bigl\{\bigl(
\sfX_i - \sfZ_i \rbH_{22}^{-1}
\rbH_{21} \bigr)' \bDelta_{i0}
\rbV_i^{-1} \bSigma_i \rbV_i^{-1}
\bDelta_{i0} \bigl(\sfX_i - \sfZ_i
\rbH_{22}^{-1}\rbH_{21} \bigr)\bigr\}
\rbH^{11},
\]
where the superscript $\vartriangle$ denotes the dependence on $\Delta_{i0}$.
%
\begin{theorem}[(Asymptotic normality)]\label{thmnormality}
The extended GEE estimator $\widehat\bbeta_V$ is asymptotically
linear, that is,
%
\begin{equation}
\label{eqgl-normality} \widehat\bbeta_V = \bbeta_0 +
\rbH^{11}\sum_{i=1}^n \bigl(
\sfX_i-\sfZ_i \rbH_{22}^{-1}
\rbH_{21}\bigr)' \bDelta_{i0}
\rbV_i^{-1} \vece_i + \RMo_P \biggl(
\frac{1}{\sqrt{n}} \biggr).
\end{equation}
Consequently,
%
\begin{equation}
\label{eqgl-asympNorm} \bigl\{\rbR^{\vartriangle}(\widehat\bbeta_V)
\bigr\}^{-1/2}(\widehat\bbeta_V - \bbeta_0)
\stackrel{d} {\longrightarrow}\operatorname{Normal}(0,\mathbf{Id}),
\end{equation}
where $\mathbf{Id}$ denotes the $K\times K$ identity matrix.
\end{theorem}

When applying the asymptotic normality result for asymptotic
inference, the variance $\rbR^{\vartriangle}(\widehat\bbeta_V)$
can be estimated by replacing $\bSigma_i$ with $(\bY_i-
\sfX_i\widehat{\bbeta}_V - \sfZ_i\widehat\bgamma) (\bY_i-
\sfX_i\widehat{\bbeta}_V - \sfZ_i\widehat\bgamma)'$, and
substituting parameter estimates in $\bDelta_{i0}$. The resulting
estimator of variance is referred to as the Sandwich estimator.
%
\begin{theorem}\label{thmgee-eff}
$\rbR^\vartriangle(\widehat\bbeta_V) \geq\rbR^\vartriangle
(\widehat
\bbeta_\Sigma)$.
\end{theorem}

Theorem \ref{thmgee-eff} says that $\widehat\bbeta_\Sigma$
is the most efficient in
the class of extended GEE estimators with general working covariance matrices.
Such a result is in parallel to that for standard parametric GEE
estimators (Liang and Zeger \cite{lz86}). This theorem is a consequence
of the
generalized Cauchy--Schwarz inequality and can be proved using exactly
the same argument as Theorem 1 of Huang, Zhang and Zhou \cite{hzz07}.

When the covariance matrices are correctly specified, the extended
GEE estimators are efficient in a stronger sense than just
described. Next, we show that the extended GEE estimator of
$\bbeta$ is the most efficient one among all regular estimators
(see Bickel\vspace*{1pt} \textit{et al.} \cite{bkrw93} for the precise definition of regular
estimators). In other words, the asymptotic variance of
$\widehat{\bbeta}_\Sigma$ achieves the semiparametric efficiency
bound, that is, the inverse of the efficient information matrix.
%
\begin{corollary}\label{corgl-eff}
The estimator $\widehat\bbeta_\Sigma$ is asymptotically normal and
semiparametric
efficient, that is,
%
\begin{equation}
\label{asyeff} (n \rbI_{\mathrm{eff}})^{1/2} (\widehat
\bbeta_\Sigma-\bbeta_0) \stackrel{d} {\longrightarrow}
\operatorname{Normal}(0,\mathbf{Id}).
\end{equation}
\end{corollary}

In the below, we sketch the proof of Corollary \ref{corgl-eff} and
postpone the details to the \hyperref[app]{Appendix}.
Fixing $\rbV_i = \bSigma_i$ in the definition of $\rbH$ as
given in (\ref{eqH}), we see that
$\rbR^{\vartriangle}(\widehat\bbeta_\Sigma)$ can be written as
\[
\rbH^{11}\sum_{i=1}^n \bigl\{
\bigl(\sfX_i - \sfZ_i \rbH_{22}^{-1}
\rbH_{21}\bigr)' \bDelta_{i0}
\bSigma_i^{-1}\bDelta_{i0} \bigl(
\sfX_i - \sfZ_i \rbH_{22}^{-1}
\rbH_{21} \bigr)\bigr\} \rbH^{11}.
\]
Using the block matrix inversion formula (\ref{eqHinv}) and
examining\vspace*{1pt} the $(1,1)$-block of the identity
$\rbH^{-1} = \rbH^{-1} \rbH\rbH^{-1}$, we obtain that
$\rbR^{\vartriangle}(\widehat\bbeta_\Sigma) = \rbH^{11}$.
Denote $\widehat\rbI_n^{-1} = n\rbR^{\vartriangle}(\widehat\bbeta_\Sigma)$.
It is easily seen using (\ref{eqHinv}) that
\begin{eqnarray*}
\widehat{\rbI}_n & = &\frac{1}{n}\sum
_{i=1}^n \sfX_i'
\bDelta_{i0} \bSigma_i^{-1} \bDelta_{i0}
\sfX_i
\\
&&{} - \frac{1}{n}\sum_{i=1}^n
\sfX_i' \bDelta_{i0} \bSigma_i^{-1}
\bDelta_{i0} \sfZ_i \Biggl(\sum
_{i=1}^n \sfZ_i'
\bDelta_{i0} \bSigma_i^{-1} \bDelta_{i0}
\sfZ_i \Biggr)^{-1}\sum_{i=1}^n
\sfZ_i' \bDelta_{i0} \bSigma_i^{-1}
\bDelta_{i0} \sfX_i.
\end{eqnarray*}
The asymptotic normality result in Theorem \ref{thmnormality}
can be rewritten as
\[
(n \widehat\rbI_n)^{1/2} (\widehat\bbeta_\Sigma-
\bbeta_0) \stackrel{d} {\longrightarrow} \operatorname{Normal} (0,
\mathbf{Id}).
\]
Thus, Corollary \ref{corgl-eff} follows from Theorem
\ref{thmnormality} and
the result that $\widehat\rbI_n \to\rbI_{\mathrm{eff}}$. The matrix
$\widehat\rbI_n$ can be interpreted as a spline-based consistent
estimate of the efficient information matrix.
%
\begin{remark}\label{remarkgl-cluster}
When $\bT$ is a cluster-level covariate, that is, $\bT_{ij}=
\widetilde
\bT_i
= (T_{i1},\ldots, T_{iD})'$,
$j=1,\ldots, m_i$, Theorems \ref{thmrate}, \ref{thmnormality}
and Corollary \ref{corgl-eff} still hold. In that case, We can
simplify C1 to the following condition.

\begin{enumerate}[C1$'$.]
\item[C1$'$.]
The random variables $T_{id}$ are bounded, uniformly in $i=1,\ldots, n$,
and $d=1,\ldots, D$.
The joint distribution of any pair of $T_{id}$ and $T_{id'}$ has a
density $f_{idd'}(t_{id},t_{id'})$ with respect to the Lebesgue measure.
We assume that $f_{idd'}(\cdot, \cdot)$ is bounded away from 0
and infinity, uniformly in $i=1,\ldots, n$, and $d,d'=1,\ldots,d$.
\end{enumerate}
\end{remark}
%
\begin{remark}\label{remarksmoothness}
Our asymptotic result on estimation of the Euclidean
parameter is quite insensitive to the choice of the
number of terms $Q_d$ in the basis expansion which plays the role
of a smoothing parameter. Specifically, suppose that additive
components of $\theta_{0,\add}(\cdot)$ and
$\varphi_{k,\add}^*(\cdot)$, $k=1,\ldots, K$, all have
bounded second derivatives, that is, condition \textup{C8$'$} is
satisfied with $\alpha=2$.
Then the requirement on $Q_d$ reduces to $n^{1/8} \ll Q_d \ll
n^{1/2}/\log^2 n$, a wide range for choosing $Q_d$.
Thus, the precise determination of $Q_d$
is not of particular concern when applying our asymptotic results.
This insensitivity of smoothing parameter
is also confirmed by our simulation study.
In practice, it is advisable to use the usual data driven methods such as
delete-cluster(subject)-out cross-validation to select $Q_d$
and then check the sensitivity of the results (Huang, Zhang and Zhou
\cite{hzz07}).
\end{remark}
%
\begin{remark}
For simplicity, we assume in our asymptotic analysis that the working
correlation parameter vector $\tau$ in $\rbV_i$ is known. It can be
estimated via the method of moments using a quadratic function of
$Y_i$'s, just as in the application of the standard parametric
GEEs (Liang and Zeger \cite{lz86}). Similar to the parametric case, as
long as such an estimate of $\tau$ converges in probability to
some $\tau^\dag$ at $\sqrt{n}$ rate, there is no asymptotic effect
on $\widehat{\bbeta}$ due to the estimation of $\tau$; see
Huang, Zhang and Zhou \cite{hzz07}, Remark 1.
\end{remark}
%
\begin{remark}
Our method does not require the assumption of normal error
distribution. However, because it is essentially a least squares
method, it is not robust to outliers. To achieve robustness to outlying
observations, it is recommended to use the M-estimator type method as
considered in He, Fung and Zhu \cite{hfz05}.
\end{remark}

\section{Numerical results}\label{sec4}

\subsection{Simulation}\label{sec4.1}
We conducted simulation studies to evaluate the finite sample
performance of the proposed method. When the number of
observations are the same per subject/cluster and the identity
link function is used, our method performs comparably to the
method of Carroll \textit{et al.} \cite{cmmy09} (see supplementary materials).
In this section, we focus on simulation setups that cannot
be handled by the existing
method of Carroll \textit{et al.} \cite{cmmy09}. We generated data from
the model
\[
E(Y_{ij}|X_{ij}, Z_{ij1}, Z_{ij2}) = g
\bigl\{\beta_0 + X_{ij}\beta_1 +
f_1(Z_{ij1}) + f_2 (Z_{ij2})\bigr\},\qquad
j=1,\ldots, n_i, i = 1,\ldots, n,
\]
where $g$ is a link function which will be specified below,
$\beta_0=0$, $\beta_1=0.5$, $f_1(t)= \sin\{2 \uppi(t-0.5)\}$,
and $f_2(t)= t - 0.5 + \sin\{2 \uppi(t-0.5)\}$.
The covariates $Z_{ij1}$ and $Z_{ij2}$ were generated from
independent Normal$(0.5,0.25)$ random variables but truncated to
the unit interval $[0,1]$. The covariate $X_{ij}$ was generated as
$X_{ij}=3(1-2Z_{ij1})(1-2Z_{ij2})+u_{ij}$ where $u_{ij}$ were
independently drawn from Normal$(0,0.25)$. We obtained different
simulation setups by varying the observational time distribution,
the correlation structure, the parameters of the correlation
function, the data distribution, and the number of subjects. We
present results for five different setups, the details of which
are given below.

For each simulation setup, 400 simulation runs were conducted and
summary statistics of the results were calculated. For each
simulated data set, the proposed generalized GEE estimator was
calculated using a working independence (WI), an exchangeable (EX)
correlation, or an autoregressive correlation structure. The
correlation parameter $\rho$ was estimated using the method of
moments. Cubic splines were used with the number of knots chosen
from the range 1--7 by the five-fold delete-subject-out
cross-validation. The bias, variance, and the mean squared errors
of Euclidean parameters were calculated for each scenario based on
the 400 runs. The mean integrated squared errors (MISE),
calculated using 100 grip points over $[0,1]$, for estimating
$f_1(\cdot)$ and $f_2(\cdot)$, were also computed.

\textit{Setup} 1. The longitudinal responses are from multivariate
normal distribution with the autoregressive correlation structure
and the identity link function. For each subject, six
observational times are evenly placed between 0 and 1. The results
are summarized in Table \ref{tbsimu1}.

\begin{table}
\caption{Summary of simulation results for setup 1, based on
400 replications. The generalized GEE estimators using a working
independence (WI), an exchangeable (EX) correlation structure, and
an autoregressive (AR) structure are compared. The true
correlation structure is the autoregressive with the lag-one
correlation being $\rho$. Each entry of the table equals the
original value multiplied by $10^5$}
\label{tbsimu1}
\begin{tabular*}{\tablewidth}{@{\extracolsep{4in minus 4in}}lld{4.0}lld{4.0}llll@{}}
\hline
& &\multicolumn{3}{l}{\hspace*{-2pt}$\beta_0=0$} &
\multicolumn{3}{l}{\hspace*{-2pt}$\beta_1=0.5$}
& \multicolumn{1}{l}{$f_1(\cdot)$} & \multicolumn{1}{l@{}}{$f_2(\cdot)$}
\\[-4pt]
& &\multicolumn{3}{l}{\hspace*{-2pt}\hrulefill} &
\multicolumn{3}{l}{\hspace*{-2pt}\hrulefill}
& \multicolumn{1}{l}{\hrulefill} & \multicolumn{1}{l@{}}{\hrulefill}
\\
$\rho$& \multicolumn{1}{l}{Method} & \multicolumn{1}{l}{\hspace*{-2pt}Bias}
& \multicolumn{1}{l}{SD} & \multicolumn{1}{l}{MSE} & \multicolumn{1}{l}{\hspace*{-2pt}Bias}
& \multicolumn{1}{l}{SD} & \multicolumn{1}{l}{MSE} & \multicolumn{1}{l}{MISE($f_1$)}
& \multicolumn{1}{l@{}}{MISE($f_2$)}%
\\
\hline
\multicolumn{10}{@{}c@{}}{$n=100$} \\
[4pt]
0.2 & WI & -27 & 391 & 391 & -176 & 60 & 60 & 1428 & 1330 \\
& EX & 14 & 381 & 381 & -173 & 58 & 58 & 1359 & 1307 \\
& AR & -27 & 373 & 373 & -180 & 57 & 57 & 1311 & 1279 \\
[4pt]
0.5 & WI & -102 & 586 & 586 & -169 & 59 & 60 & 1452 & 1322 \\
& EX & 44 & 524 & 524 & -168 & 48 & 49 & 1245 & 1116 \\
& AR & -42 & 474 & 474 & -152 & 40 & 40 & \hphantom{0}990 & \hphantom{0}920 \\
[4pt]
0.8 & WI & -194 & 924 & 925 & -100 & 60 & 60 & 1448 & 1358 \\
& EX & -13 & 787 & 787 & -133 & 29 & 29 & \hphantom{0}747 & \hphantom{0}662 \\
& AR & -96 & 686 & 686 & -93 & 16 & 16 & \hphantom{0}463 & \hphantom{0}461 \\
[4pt]
\multicolumn{10}{@{}c@{}}{$n=200$} \\
[4pt]
0.2 & WI & -239 & 181 & 182 & -35 & 26 & 26 & \hphantom{0}689 & \hphantom{0}709 \\
& EX & -273 & 180 & 181 & -47 & 25 & 25 & \hphantom{0}669 & \hphantom{0}698 \\
& AR & -238 & 175 & 176 & -32 & 26 & 26 & \hphantom{0}656 & \hphantom{0}664 \\
[4pt]
0.5 & WI & -261 & 270 & 271 & 4 & 27 & 27 & \hphantom{0}676 & \hphantom{0}712 \\
& EX & -281 & 258 & 259 & -38 & 23 & 23 & \hphantom{0}569 & \hphantom{0}604 \\
& AR & -208 & 241 & 242 & -18 & 19 & 19 & \hphantom{0}482 & \hphantom{0}493 \\
[4pt]
0.8 & WI & -183 & 448 & 449 & 62 & 30 & 30 & \hphantom{0}677 & \hphantom{0}723 \\
& EX & -243 & 400 & 401 & -20 & 13 & 13 & \hphantom{0}338 & \hphantom{0}361 \\
& AR & -162 & 369 & 369 & -7 & 8 & 8 & \hphantom{0}224 & \hphantom{0}245 \\
\hline
\end{tabular*}
\end{table}

\textit{Setup} 2. The same as setup 1, except that the log link is
used. The results are summarized in Table \ref{tbsimu2}.

\begin{table}
\caption{Summary of simulation results for setup 2, based on
400 replications. The generalized GEE estimators using a working
independence (WI), an exchangeable (EX) correlation structure, and
an autoregressive (AR) structure are compared. The true
correlation structure is the autoregressive with the lag-one
correlation being $\rho$. Each entry of the table equals the
original value multiplied by $10^5$}
\label{tbsimu2}
\begin{tabular*}{\tablewidth}{@{\extracolsep{4in minus 4in}}lld{5.0}d{4.0}rd{3.0}lld{4.0}d{4.0}@{}}
\hline
& &\multicolumn{3}{l}{\hspace*{-2pt}$\beta_0=0$} &
\multicolumn{3}{l}{$\beta_1=0.5$}
& \multicolumn{1}{l}{$f_1(\cdot)$} & \multicolumn{1}{l@{}}{$f_2(\cdot)$}
\\[-4pt]
& &\multicolumn{3}{l}{\hspace*{-2pt}\hrulefill} &
\multicolumn{3}{l}{\hrulefill}
& \multicolumn{1}{l}{\hrulefill} & \multicolumn{1}{l@{}}{\hrulefill}
\\
$\rho$& \multicolumn{1}{l}{Method} & \multicolumn{1}{l}{\hspace*{-2pt}Bias}
& \multicolumn{1}{l}{SD} & \multicolumn{1}{l}{MSE} & \multicolumn{1}{l}{Bias}
& \multicolumn{1}{l}{SD} & \multicolumn{1}{l}{MSE} & \multicolumn{1}{l}{MISE($f_1$)}
& \multicolumn{1}{l@{}}{MISE($f_2$)}%
\\
\hline
\multicolumn{10}{@{}c@{}}{$n=100$}\\
[4pt]
0.2 & WI & -2294 & 970 & 1022 & 44 & 30 & 30 & 1407 & 2280 \\
& EX & -2223 & 962 & 1011 & 77 & 29 & 29 & 1397 & 2195 \\
& AR & -2265 & 951 & 1003 & 64 & 29 & 29 & 1374 & 2104 \\
[4pt]
0.5 & WI & -2164 & 1137 & 1183 & 62 & 33 & 33 & 1356 & 2198 \\
& EX & -1928 & 1007 & 1045 & 117 & 26 & 26 & 1146 & 1846 \\
& AR & -1711 & 799 & 828 & 89 & 23 & 23 & 948 & 1394 \\
[4pt]
0.8 & WI & -2361 & 1727 & 1783 & 85 & 37 & 37 & 1378 & 2234 \\
& EX & -2091 & 1017 & 1061 & 140 & 17 & 17 & 742 & 1303 \\
& AR & -1911 & 722 & 758 & 116 & 12 & 12 & 518 & 824 \\
[4pt]
\multicolumn{10}{@{}c@{}}{$n=200$} \\
[4pt]
0.2 & WI & -1387 & 388 & 407 & 88 & 16 & 16 & 601 & 1010 \\
& EX & -1498 & 390 & 412 & 99 & 16 & 16 & 582 & 1024 \\
& AR & -1497 & 384 & 407 & 98 & 15 & 15 & 565 & 985 \\
[4pt]
0.5 & WI & -1433 & 499 & 519 & 84 & 17 & 17 & 618 & 1068 \\
& EX & -1532 & 435 & 458 & 108 & 14 & 14 & 534 & 830 \\
& AR & -1525 & 387 & 410 & 97 & 12 & 12 & 436 & 660 \\
[4pt]
0.8 & WI & -1410 & 712 & 732 & 76 & 18 & 18 & 623 & 1095 \\
& EX & -1192 & 332 & 346 & 81 & \hphantom{0}9 & \hphantom{0}9 & 302 & 482 \\
& AR & -1433 & 277 & 298 & 88 & \hphantom{0}5 & \hphantom{0}5 & 219 & 323 \\
\hline
\end{tabular*}
\end{table}

\textit{Setup} 3. This setup is the same as setup 1, except that the
exchangeable correlation structure is used and the observational
time distribution is different. For each subject, ten
observational times are first evenly placed between 0 and 1. Then
40\% of the observations are removed from each dataset and thus
different subjects may have different number of observations and
the observational times may be irregularly placed. The results are
summarized in Table \ref{tbsimu3}.

\begin{table}
\caption{Summary of simulation results for setup 3, based on 400
replications. The generalized GEE estimators using a working
independence (WI) and an exchangeable (EX) correlation structure
are compared. The true correlation structure is the exchangeable
with parameter $\rho$. Each entry of the table equals the
original value multiplied by $10^5$}
\label{tbsimu3}
\begin{tabular*}{\tablewidth}{@{\extracolsep{4in minus 4in}}lld{4.0}d{4.0}rd{4.0}lld{4.0}d{4.0}@{}}
\hline
& &\multicolumn{3}{l}{\hspace*{-2pt}$\beta_0=0$} &
\multicolumn{3}{l}{\hspace*{-2pt}$\beta_1=0.5$}
& \multicolumn{1}{l}{$f_1(\cdot)$} & \multicolumn{1}{l@{}}{$f_2(\cdot)$}
\\[-4pt]
& &\multicolumn{3}{l}{\hspace*{-2pt}\hrulefill} &
\multicolumn{3}{l}{\hspace*{-2pt}\hrulefill}
& \multicolumn{1}{l}{\hrulefill} & \multicolumn{1}{l@{}}{\hrulefill}
\\
$\rho$& \multicolumn{1}{l}{Method} & \multicolumn{1}{l}{\hspace*{-2pt}Bias}
& \multicolumn{1}{l}{SD} & \multicolumn{1}{l}{MSE} & \multicolumn{1}{l}{\hspace*{-2pt}Bias}
& \multicolumn{1}{l}{SD} & \multicolumn{1}{l}{MSE} & \multicolumn{1}{l}{MISE($f_1$)}
& \multicolumn{1}{l@{}}{MISE($f_2$)}%
\\
\hline
\multicolumn{10}{@{}c@{}}{$n=100$}\\
[4pt]
0 & WI & -129 & 337 & 337 & -125 & 61 & 61 & 1426 & 1412 \\
& EX & -109 & 336 & 336 & -116 & 61 & 61 & 1416 & 1410 \\
[4pt]
0.2 & WI & -96 & 527 & 527 & -56 & 61 & 61 & 1445 & 1423 \\
& EX & -161 & 511 & 511 & -66 & 55 & 55 & 1297 & 1347 \\
[4pt]
0.5 & WI & -125 & 797 & 798 & 14 & 62 & 62 & 1515 & 1399 \\
& EX & -216 & 735 & 735 & -39 & 37 & 37 & 924 & 962 \\
[4pt]
0.8 & WI & -23 & 1054 & 1054 & 58 & 62 & 62 & 1552 & 1362 \\
& EX & -164 & 914 & 914 & -27 & 15 & 15 & 455 & 464 \\
[4pt]
\multicolumn{10}{@{}c@{}}{$n=200$} \\
[4pt]
0 & WI & 48 & 149 & 149 & 70 & 29 & 29 & 780 & 649 \\
& EX & 54 & 149 & 149 & 74 & 29 & 29 & 782 & 659 \\
[4pt]
0.2 & WI & -99 & 253 & 253 & 39 & 29 & 29 & 798 & 677 \\
& EX & -21 & 237 & 237 & 37 & 25 & 25 & 693 & 609 \\
[4pt]
0.5 & WI & -192 & 403 & 404 & -3 & 31 & 31 & 768 & 690 \\
& EX & -64 & 354 & 354 & 15 & 16 & 16 & 470 & 432 \\
[4pt]
0.8 & WI & -240 & 564 & 565 & -60 & 32 & 32 & 718 & 702 \\
& EX & -96 & 466 & 466 & 6 & \hphantom{0}7 & \hphantom{0}7 & 236 & 227 \\
\hline
\end{tabular*}\vspace*{-3pt}
\end{table}

\textit{Setup} 4. It is the same as setup 3, except that the log
link is used. The results are summarized in Table \ref{tbsimu4}.

\begin{table}
\caption{
Summary of simulation results for setup 4, based on 400
replications. The generalized GEE estimators using a working
independence (WI) and an exchangeable (EX) correlation structure
are compared. The true correlation structure is the exchangeable
with parameter $\rho$. Each entry of the table equals the
original value multiplied by $10^5$} \label{tbsimu4}
\begin{tabular*}{\tablewidth}{@{\extracolsep{4in minus 4in}}lld{5.0}d{4.0}rd{4.0}d{2.0}d{2.0}d{4.0}d{4.0}@{}}
\hline
& &\multicolumn{3}{l}{\hspace*{-2pt}$\beta_0=0$} &
\multicolumn{3}{l}{\hspace*{-2pt}$\beta_1=0.5$}
& \multicolumn{1}{l}{$f_1(\cdot)$} & \multicolumn{1}{l@{}}{$f_2(\cdot)$}
\\[-4pt]
& &\multicolumn{3}{l}{\hspace*{-2pt}\hrulefill} &
\multicolumn{3}{l}{\hspace*{-2pt}\hrulefill}
& \multicolumn{1}{l}{\hrulefill} & \multicolumn{1}{l@{}}{\hrulefill}
\\
$\rho$& \multicolumn{1}{l}{Method} & \multicolumn{1}{l}{\hspace*{-2pt}Bias}
& \multicolumn{1}{l}{SD} & \multicolumn{1}{l}{MSE} & \multicolumn{1}{l}{\hspace*{-2pt}Bias}
& \multicolumn{1}{l}{SD} & \multicolumn{1}{l}{MSE} & \multicolumn{1}{l}{MISE($f_1$)}
& \multicolumn{1}{l@{}}{MISE($f_2$)}%
\\
\hline
\multicolumn{10}{@{}c@{}}{$n=100$}\\
[4pt]
0 & WI & -2451 & 756 & 816 & -218 & 27 & 28 & 1313 & 2318 \\
& EX & -2500 & 773 & 835 & -229 & 27 & 28 & 1329 & 2314 \\
[4pt]
0.2 & WI & -2581 & 1054 & 1120 & -176 & 30 & 31 & 1461 & 2184 \\
& EX & -2440 & 974 & 1034 & -164 & 26 & 26 & 1240 & 2012 \\
[4pt]
0.5 & WI & -2455 & 1514 & 1574 & -88 & 36 & 36 & 1574 & 2287 \\
& EX & -1970 & 918 & 956 & -102 & 19 & 19 & 851 & 1482 \\
[4pt]
0.8 & WI & -2520 & 2029 & 2093 & 3 & 41 & 41 & 1806 & 2342 \\
& EX & -2547 & 1036 & 1101 & 24 & 10 & 10 & 629 & 1025 \\
[4pt]
\multicolumn{10}{@{}c@{}}{$n=200$}\\
[4pt]
0 & WI & -883 & 329 & 336 & 114 & 14 & 14 & 653 & 826 \\
& EX & -866 & 329 & 337 & 116 & 14 & 14 & 655 & 823 \\
[4pt]
0.2 & WI & -1090 & 475 & 487 & 84 & 14 & 14 & 728 & 844 \\
& EX & -903 & 390 & 398 & 84 & 13 & 13 & 631 & 734 \\
[4pt]
0.5 & WI & -1310 & 718 & 736 & 44 & 16 & 16 & 779 & 932 \\
& EX & -951 & 344 & 353 & 67 & 9 & 9 & 421 & 538 \\
[4pt]
0.8 & WI & -1533 & 966 & 989 & 13 & 18 & 18 & 744 & 1086 \\
& EX & -1245 & 285 & 301 & 65 & 4 & 4 & 209 & 381 \\
\hline
\end{tabular*}
\end{table}

\textit{Setup} 5. This setup is the same as setup 4, except that the
Poisson distribution is used as the marginal distribution.
All regression parameters in the general setup, the Euclidean and
the functional, are halved for appropriate scaling of the response
variable. The results are summarized in Table \ref{tbsimu5}.

\begin{table}[t!]
\caption{Summary of simulation results for setup 5, based on 400
replications. The generalized GEE estimators using a working
independence (WI) and an exchangeable (EX) correlation structure
are compared. The true correlation structure is the exchangeable
with parameter $\rho$. Each entry of the table equals the
original value multiplied by $10^5$}
\label{tbsimu5}
\begin{tabular*}{\tablewidth}{@{\extracolsep{4in minus 4in}}lld{5.0}d{4.0}rd{4.0}d{3.0}
d{3.0}d{4.0}d{4.0}@{}}
\hline
& &\multicolumn{3}{l}{\hspace*{-2pt}$\beta_0=0$} &
\multicolumn{3}{l}{\hspace*{-2pt}$\beta_1=0.5$}
& \multicolumn{1}{l}{$f_1(\cdot)$} & \multicolumn{1}{l@{}}{$f_2(\cdot)$}
\\[-4pt]
& &\multicolumn{3}{l}{\hspace*{-2pt}\hrulefill} &
\multicolumn{3}{l}{\hspace*{-2pt}\hrulefill}
& \multicolumn{1}{l}{\hrulefill} & \multicolumn{1}{l@{}}{\hrulefill}
\\
$\rho$& \multicolumn{1}{l}{Method} & \multicolumn{1}{l}{\hspace*{-2pt}Bias}
& \multicolumn{1}{l}{SD} & \multicolumn{1}{l}{MSE} & \multicolumn{1}{l}{\hspace*{-2pt}Bias}
& \multicolumn{1}{l}{SD} & \multicolumn{1}{l}{MSE} & \multicolumn{1}{l}{MISE($f_1$)}
& \multicolumn{1}{l@{}}{MISE($f_2$)}%
\\
\hline
\multicolumn{10}{@{}c@{}}{$n=100$}\\
[4pt]
0 & WI & -2967 & 366 & 454 & -379 & 113 & 114 & 1376 & 1446 \\
& EX & -2983 & 368 & 457 & -353 & 112 & 113 & 1368 & 1439 \\
[4pt]
0.2 & WI & -2557 & 738 & 803 & -456 & 120 & 122 & 1394 & 1385 \\
& EX & -2998 & 777 & 867 & -367 & 98 & 99 & 1031 & 1110 \\
[4pt]
0.5 & WI & -1952 & 1101 & 1140 & 221 & 126 & 126 & 1446 & 1484 \\
& EX & -2272 & 1339 & 1390 & 215 & 70 & 70 & 506 & 628 \\
[4pt]
0.8 & WI & -1979 & 1344 & 1383 & 369 & 126 & 127 & 1464 & 1567 \\
& EX & -2349 & 1651 & 1706 & 506 & 71 & 74 & 411 & 545 \\
[4pt]
\multicolumn{10}{@{}c@{}}{$n=200$}\\
[4pt]
0 & WI & -1563 & 190 & 214 & -214 & 51 & 52 & 685 & 780 \\
& EX & -1586 & 191 & 216 & -208 & 51 & 52 & 691 & 781 \\
[4pt]
0.2 & WI & -1015 & 405 & 415 & -195 & 54 & 55 & 637 & 771 \\
& EX & -1355 & 402 & 421 & -154 & 45 & 45 & 516 & 589 \\
[4pt]
0.5 & WI & -1301 & 599 & 616 & 143 & 55 & 55 & 742 & 777 \\
& EX & -1751 & 634 & 665 & 218 & 30 & 30 & 256 & 300 \\
[4pt]
0.8 & WI & -1381 & 636 & 655 & 341 & 52 & 53 & 768 & 802 \\
& EX & -1942 & 662 & 699 & 434 & 30 & 32 & 224 & 281 \\
\hline
\end{tabular*} \vspace*{6pt}
\end{table}

\begin{table}[t!]
\caption{Estimates of the Euclidean parameters in the CD4 cell
counts study using the spline-based estimates. Working correlation
structures used are working independence (WI) and exchangeable
(EX). The standard errors (SE) are calculated using the sandwich
formula} \label{tbcd4}
\begin{tabular*}{\tablewidth}{@{\extracolsep{\fill}}lllll@{}}
\hline
&\multicolumn{2}{l}{WI} &\multicolumn{2}{l@{}}{EX}\\
[-4pt]
&\multicolumn{2}{l}{\hrulefill} &\multicolumn{2}{l@{}}{\hrulefill}\\
Parameter & \multicolumn{1}{l}{Estimate} & SE & \multicolumn{1}{l}{Estimate}
& SE \\
\hline
Smoking &\hphantom{$-$}0.0786 &0.0119 &\hphantom{$-$}0.0619 &0.0111\\
Drug &\hphantom{$-$}0.0485&0.0421 &\hphantom{$-$}0.0134 &0.0294\\
Sex partners & $-$0.0056 &0.0043 &\hphantom{$-$}0.0017 &0.0035\\
Depression & $-$0.0025 &0.0014 & $-$0.0031 &0.0013\\
\hline
\end{tabular*}
\end{table}

We have the following observations from the simulation results:
for both Euclidean parameters, the estimator accounting for the
correlation is more efficient (and sometimes significantly so)
than the estimator using working independence correlation
structure, even when the correlation structure is misspecified.
Using the correct correlation structure usually produces the most
efficient estimation. Efficiency gain gets bigger when the
correlation parameter $\rho$ gets larger. The variance is usually
a dominating factor when comparing the MSEs between the two
estimators. We have also observed that the sandwich estimated SEs
work reasonably well; the averages of the sandwich estimated SEs
are close to the Monte Carlo sample standard deviations (numbers
not shown to save space). For the functional parameters
$f_1(\cdot)$ and $f_2(\cdot)$, the spline estimator accounting for
the correlation is more efficient and the most efficient when the
working correlation is the same as the true correlation structure.
We also examined the Normal Q--Q plots of the Euclidean parameter
estimates and observed that the distributions of the estimates are
close to normal. These empirical results agree nicely with our
theoretical results.

\subsection{The longitudinal CD4 cell count data}\label{sec4.2}

To illustrate our method on a real data set, we considered the
longitudinal CD4 cell count data among HIV seroconverters
previously analyzed by Zeger and Diggle \cite{ZeDi94}. This data set
contains 2376 observations of CD4$+$ cell counts on 369 men infected
with the HIV virus. See Zeger and Diggle \cite{ZeDi94} for more detailed
description of the data. We fit a partially linear additive model
using the log link with the CD4 counts as the response, covariates
entering the model linearly including smoking status measured by
packs of cigarettes, drug use (yes, 1; no 0), number of sex
partners, and depression status measures by the CESD scale (large
values indicating more depression symptoms), and the effects of
age and time since seroconversion being modeled nonparametrically.
We would like to remark that the partially linear additive
model here provides a good balance of model interpretability and flexibility.
Age and time are of continuous type and thus their effects are naturally
modeled nonparametrically. Other variables are of discrete type
and are not suitable for a nonparametric model.

Table \ref{tbcd4} gives the estimates of the Euclidean parameters
using both the WI and EX correlation structures. Cubic splines
were used for fitting the additive functions and reported results
correspond to the number of knots selected by the five-fold
delete-subject-out cross-validation from the range of 0--10. The
selected numbers of knots are 8 for time and 4 for age when using
the WI structure and 8 for time and 3 for age when using the EX
structure. The estimates of the Euclidean parameters using the EX
structure have smaller SE than those using the WI structure,
suggesting that the EX structure produces more efficient estimates
for this data set.

\begin{appendix}\label{app}
\section*{Appendix}
\setcounter{subsection}{0}
\subsection{\texorpdfstring{Proof of Lemma \protect\ref{lemmaeff-score}
(derivation of the efficient score)}
{Proof of Lemma 1 (derivation of the efficient score)}}\label{secA.1}
Let $\dot\ell_{\beta}$ denote\vspace*{2pt} the ordinary score for $\bbeta$ when only
$\bbeta$ is unknown.
Let $\mathcal{P}_{f}$ and
$\mathcal{P}_{\theta}$ be the models with only
$\{f_i, i=1,\ldots,n\}$ and $\theta_{\add}(\cdot)$ unknown,
respectively, and let $\dot{\mathcal{P}}_{f}$ and $\dot{\mathcal
{P}}_{\theta}$
be the corresponding tangent spaces.\vadjust{\goodbreak}
Following the discussions in Section 3.4 of Bickel \textit{et al.} \cite{bkrw93}
(see also Appendix A6 of Huang, Zhang and Zhou \cite{hzz07}), we have
%
\renewcommand{\theequation}{\arabic{equation}}
\begin{equation}
\label{effscocal} \ell^\ast_\beta= \Pi\bigl[\dot
\ell_\beta|\dot{\mathcal{P}}_{f}^\bot\bigr] -
\Pi\bigl[\Pi\bigl(\dot\ell_\beta|\dot{\mathcal{P}}_{f}^\bot
\bigr) | \Pi\bigl[\dot{\mathcal{P}}_\theta|\dot{\mathcal{P}}_{f}^\bot
\bigr] \bigr],
\end{equation}
where $\Pi[\cdot|\cdot]$ denote the projection operator,
and $\dot{\mathcal{P}}^\bot$ denote the orthogonal complement of
$\dot{\mathcal{P}}$.
Lemma A.4 in Huang, Zhang and Zhou \cite{hzz07} directly implies that
%
\begin{equation}
\label{effscocal1} \Pi\bigl[\dot\ell_\beta|\dot{\mathcal{P}}_{f}^\bot
\bigr] = \sum_{i=1}^n {
\sfX}_i'\bDelta_{i0}\bSigma_i^{-1}
\bigl[{\bY}_i- \mu\bigl({\sfX}_i\bbeta_0+
\theta_{0,\add}({\sfT}_i)\bigr)\bigr].
\end{equation}
Similarly, by constructing parametric submodels for each
$\theta_k(\cdot)$ and slightly adapting the same Lemma, we have
%
\begin{equation}
\label{effscocal2} \Pi\bigl[\dot{\mathcal{P}}_\theta|\dot{\mathcal
{P}}_{f}^\bot\bigr] = \sum_{i=1}^{n}
\Biggl(\sum_{d=1}^{D} \psi_d({
\bT}_{id}) \Biggr)' \bDelta_{i0}
\bSigma_i^{-1}\bigl[{\bY}_i-\mu\bigl({\sfX
}_i\bbeta_0 +\theta_{0,\add}({
\sfT}_i)\bigr)\bigr],
\end{equation}
where $\psi_d(\bT_{id})=(\psi_d(T_{i1d}),\ldots,\psi
_{im_id}(T_{im_id}))'$, for $\psi_d(\cdot)\in L_2(\mathcal{T}_d)$.
Combination of (\ref{effscocal})--(\ref{effscocal2}) gives (\ref
{eqeffsco}).

\subsection{\texorpdfstring{Proof sketch for Theorem \protect\ref{thmrate} (consistency)}
{Proof sketch for Theorem 1 (consistency)}}\label{secA.2}
Let $\epsilon_n=(Q_n/n)^{1/2}\log n\vee\rho_n$. To show
(\ref{convrate}), it suffices to show that $P(\|\widehat
f_n-f_n^\ast\|_n>\epsilon_n)\rightarrow0$ as
$n\rightarrow\infty$. Applying the peeling device
(see the proof of Theorem 9.1 of van de Geer \cite{vg00}), we can bound
the above
probability by the sum of $2C_0\exp(-n\epsilon_n^2/(256C_0^2))$
and $P(\|y-f_n^\ast\|_n>\sigma)$ for some positive constant $C_0$.
Considering condition C8 and choosing some proper $\sigma$ related
to~$\rho_n$, we complete the proof of (\ref{convrate}). As for
(\ref{noncons}), we have that
\[
\bigl\|\widehat f_n - f_n^*\bigr\|_\infty\aplt\bigl\|
\mathbf{x}'\widehat{\bbeta}_V + \widehat\theta- \bigl(
\mathbf{x}' \bbeta_0 + \theta^*_n\bigr)
\bigr\|_\infty\aplt Q_n^{1/2} \bigl\|\mathbf{x}'
\widehat{\bbeta}_V + \widehat\theta- \bigl(\mathbf{x}'
\bbeta_0 + \theta^*_n\bigr)\bigr\|
\]
by Condition C5(iii) and Lemma S.2 in the supplementary note that
%
\begin{equation}\label{inter2}
\bigl\|\mathbf{x}'\beta+ g(\mathbf{t})\bigr\|_{\infty}
\aplt{Q_n}^{1/2} \bigl\| \mathbf{x}'\beta+ g(
\mathbf{t})\bigr\| \qquad\mbox{for } g\in\mathbb{G}_\add.
\end{equation}
It then follows by condition C5(ii) and (\ref{convrate}) that
\[
Q_n^{1/2} \bigl\|\mathbf{x}'\widehat{
\bbeta}_V + \widehat\theta- \bigl(\mathbf{x}'
\bbeta_0 + \theta^*_n\bigr)\bigr\| \aplt Q_n^{1/2}
\RMO_P\bigl\{(Q_n/n)^{1/2}\log n +
\rho_n\bigr\} = \RMo_P(1)
\]
since $(Q_n\log n)^2/n\rightarrow0$ and $Q_n\rho_n^2\rightarrow
0$ by condition C8 and the fact that $\rho_n\asymp Q_n^{-\alpha}$
for $\alpha>1/2$. We thus obtain (\ref{noncons}). Due to condition
C5(iii), it follows that $\|f_n^* - f_0\|_\infty= \RMO(\|\theta_n^*
- \theta_{0,\add}\|_\infty) = \RMO(\rho_n)$ by Taylor's theorem.
Combining this with (\ref{noncons}), we obtain
(\ref{eqconsis-whole}). From the proof of (\ref{noncons}), we
have that
\[
\bigl\|\mathbf{x}'\widehat{\bbeta}_V + \widehat\theta-
\bigl(\mathbf{x}' \bbeta_0 + \theta^*_n
\bigr)\bigr\| = \RMO_P\bigl\{(Q_n/n)^{1/2}\log n +
\rho_n\bigr\}.
\]
Considering Lemma 3.1 of Stone \cite{s94}, we obtain that
$\|\mathbf{x}'(\widehat{\bbeta}_V- \bbeta_0)\|^2 = \RMo_P(1)$, which
together with the no-multicollinearity condition C2 implies
$\widehat\bbeta_V \stackrel{P}{\to} \bbeta_0$. By (\ref{inter2}),
we also obtain
\[
\bigl\|\widehat\theta- \theta^*_n\bigr\|_\infty= Q_n^{1/2}
\RMO_P\bigl\{(Q_n/n)^{1/2}\log n +
\rho_n\bigr\} = \RMo_P(1).
\]
Since $\|\theta_n^* - \theta_{0,\add}\|_\infty= \RMO(\rho_n) =
\RMo(1)$, application of the triangle inequality yields
$\|\widehat\theta- \theta_{0,\add}\|_\infty= \RMo_P(1)$, the last
conclusion.

\subsection{\texorpdfstring{Proof sketch for Theorem \protect\ref{thmnormality} (asymptotic normality)}
{Proof sketch for Theorem 2 (asymptotic normality)}}\label{secA.3}
Note that $\widehat\bbeta_V \in\bbR^{K}$ and $\widehat
\bgamma\in\bbR^{Q_n}$ solve the estimating equations
%
\begin{equation}
\label{eqgl-estimatingEquation-0} \sum_{i=1}^n
\sfU_i' \widehat\bDelta_i
\rbV_i^{-1} \bigl\{\bY_i - \mu(
\sfX_i\widehat{\bbeta}_V + \sfZ_i\widehat
\bgamma)\bigr\} = 0
\end{equation}
with $\sfU_i=(\sfX_i,\sfZ_i)$, and $\widehat\bDelta_i$ is a
diagonal matrix with the diagonal elements being the first
derivative of $\mu(\cdot)$ evaluated at $X_{ij}'\widehat{\bbeta}_V
+ Z_{ij}'\widehat\bgamma$, $j=1,\ldots, m_i$. Using the Taylor
expansion, we have that
%
\begin{equation}
\label{eqgl-taylor-0} \mu(\sfX_i\widehat{\bbeta}_V +
\sfZ_i\widehat\bgamma)\approx\mu\bigl(\sfX_i
\bbeta_0 + \theta_0(\bT_i) \bigr) +
\bDelta_{i0} \bigl\{\sfX_i(\widehat{\bbeta}_V-
\bbeta_0)+ \sfZ_i\widehat\bgamma-
\theta_0(\bT_i) \bigr\}.
\end{equation}
Recall that $\bgamma^\ast$ is assumed to satisfy
$\rho_n=\|\theta_{0,\add}- \bB'\bgamma^\ast\|_{\infty
}\rightarrow
0$. Substituting (\ref{eqgl-taylor-0}) into
(\ref{eqgl-estimatingEquation-0}) yields
%
\begin{equation}
\label{eqgl-estimatingEquation1-0} 0 = \sum_{i=1}^n
\sfU_i' (\widetilde\bJ_1 + \widetilde
\bJ_2) - \sum_{i=1}^n
\sfU_i' \bDelta_{i0} \rbV_i^{-1}
\bDelta_{i0} \sfU_i %
\pmatrix{\widehat{
\bbeta}_V-\bbeta_0
\cr
\widehat\bgamma-\bgamma^*},
\end{equation}
where
\[
\widetilde\bJ_1 = (\widehat\bDelta_i -
\bDelta_{i0}) \rbV_i^{-1} \bigl\{
\bY_i - \mu(\sfX_i\widehat{\bbeta}_V +
\sfZ_i\widehat\bgamma) \bigr\}
\]
and
\[
\widetilde\bJ_2 = \bDelta_{i0} \rbV_i^{-1}
\bigl\{\bY_i - \mu\bigl(\sfX_i\bbeta_0 +
\theta_0(\bT_i) \bigr) - \bDelta_{i0} \bigl(
\sfZ_i \bgamma^* - \theta_0(\bT_i) \bigr)
\bigr\}.
\]
Recalling (\ref{eqH}) and using (\ref{eqHinv}), we obtain from
(\ref{eqgl-estimatingEquation1-0}) that
\begin{eqnarray*}
\label{eqasymp-linear} \widehat{\bbeta}_V & = &\bbeta_0
+ \rbH^{11} \sum_{i=1}^n \bigl(
\sfX_i-\sfZ_i \rbH_{22}^{-1}
\rbH_{21}\bigr)' (\widetilde\bJ_1 +\widetilde
\bJ_2)
\\
&=& \bbeta_0 + \rbH^{11}\sum
_{i=1}^n \bigl(\sfX_i-
\sfZ_i \rbH_{22}^{-1} \rbH_{21}
\bigr)' \bDelta_{i0} \rbV_i^{-1}
\vece_i + \pi_n,
\end{eqnarray*}
where the error term $\pi_n$ has an explicit form and can be shown to
be $\RMo_P(n^{-1/2})$ (the proof of this part relies heavily on the
empirical process
theory and is very lengthy). By the asymptotic linear expansion
(\ref{eqgl-normality}), we have
\begin{eqnarray*}
&& \bigl\{\rbR^\vartriangle(\widehat{\bbeta}_V)\bigr
\}^{-1/2}(\widehat{\bbeta}_V-\bbeta_0)
\\
&&\quad =\bigl\{\rbR^\vartriangle(\widehat{\bbeta}_V)\bigr
\}^{-1/2}\Biggl(\rbH^{11}\sum_{i=1}^n
\bigl(\sfX_i-\sfZ_i \rbH_{22}^{-1}
\rbH_{21}\bigr)' \bDelta_{i0}
\rbV_i^{-1} \vece_i\Biggr)+
\RMo_P(1).
\end{eqnarray*}
Then by applying the central limit theorem to the above equation
and using the fact that
\[
\var\Biggl(\rbH^{11}\sum_{i=1}^n
\bigl(\sfX_i-\sfZ_i \rbH_{22}^{-1}
\rbH_{21}\bigr)' \bDelta_{i0}
\rbV_i^{-1} \vece_i \Big|\{\sfX_i,
\sfT_i\}_{i=1}^n \Biggr) = \rbR^{\vartriangle}(
\widehat\bbeta_V),
\]
we complete the whole proof of (\ref{eqgl-asympNorm}).

\subsection{\texorpdfstring{Proof of Corollary \protect\ref{corgl-eff}}
{Proof of Corollary 1}}\label{secA.4}\vspace*{-2pt}
We only need to show that $\widehat\rbI_n \to\rbI_{\mathrm{eff}}$.
Fix $\rbV_i=\bSigma_i$ in the definitions of
$\langle\xi_1,\xi_2\rangle_n^\vartriangle$ and
$\langle\xi_1,\xi_2\rangle^\vartriangle$.
Let $\widehat{\psi}_{k,n}=\arg\min_{\psi\in\mathbb{G}_\add}\|x_k
-\psi\|_{n}^{\vartriangle}$.
From (\ref{eqHinv}), we see that $\widehat\rbI_n
= (\rbH_{11} - \rbH_{12}\rbH_{22}^{-1}\rbH_{21})/n$.
Thus,
the $(k,k')$th element of $\widehat{\rbI}_n$ is
$\langle x_k-\widehat{\psi}_{k,n}, x_{k'}
-\widehat{\psi}_{k',n}\rangle_{n}^\vartriangle$.
On the other hand, by (\ref{psikast}) and~(\ref{eqinfo}),
the $(k,k')$th element of ${\rbI}_{\mathrm{eff}}$
is the limit of
$\langle x_k- {\psi}_{k}^*, x_{k'}
- {\psi}_{k'}^* \rangle_n^\vartriangle$,
where $\psi_k^* = \psi_{k,\add}^* = \argmin_{L_{2,\add}} \|x_k -
\psi\|^{\vartriangle}$.
Hence, it suffices to show that
%
\begin{equation}
\label{eqapproxI-0} \bigl\|\widehat{\psi}_{k,n}-\psi_k^\ast
\bigr\|_{n}^\vartriangle=\RMo_{P}(1),\qquad k=1,2,\ldots,K,
\end{equation}
because, if this is true, then by the triangle
inequality,
\begin{eqnarray*}
\widehat{\rbI}_n\bigl(k,k'\bigr) & = & \langle
x_k-\widehat{\psi}_{k,n}, x_{k'} -\widehat{
\psi}_{k',n}\rangle_{n}^\vartriangle
\\
& = &\bigl\langle x_k-\psi_k^\ast,x_{k'}-
\psi^\ast_{k'}\bigr\rangle_{n}^\vartriangle+\RMo_{P}(1)
= \rbI_{\mathrm{eff}}\bigl(k,k'\bigr)+\RMo_{P}(1).
\end{eqnarray*}

To show (\ref{eqapproxI-0}), we use $\psi_{k,n}^* = \Pi_n^\vartriangle x_k$
as a bridge. Notice that
\[
\bigl\|\widehat\psi_{k,n}-\psi_k^*\bigr\|_n^\vartriangle
\leq\bigl\|\psi_{k,n}^*-\psi_k^*\bigr\|_n^\vartriangle+
\bigl\|\widehat\psi_{k,n}-\psi_{k,n}^*\bigr\|_n^\vartriangle.
\]
We inspect separately the sizes of the two terms on the right-hand side
of the above inequality. First note that $\psi_{k,n}^* =
\Pi_n^\vartriangle\psi_k^*$ since $\mathbb{G}_\add\subset
L_{2,\add}$.
Thus, $\|\psi_{k,n}^*-\psi_k^*\|^\vartriangle=
\inf_{g\in\bbG_\add}\|g-\psi_k^*\|^\vartriangle\asymp
\inf_{g\in\bbG_\add}\|g-\psi_k^*\|_{L_2} = \RMO(\rho_n) =\RMo(1)$,
using Lemma S.2 in the supplementary note. Since
$E(\{\|\psi_{k,n}^*-\psi_k^*\|_n^\vartriangle\}^2) =
\{\|\psi_{k,n}^*-\psi_k^*\|^\vartriangle\}^2$, we have that
$\|\psi_{k,n}^*-\psi_k^*\|_n^\vartriangle= \RMo_P(1)$. On
the other hand, since $\psi_{k,n}^* = \Pi_n^\vartriangle x_k$ and
$\widehat\psi_{k,n} =\widehat\Pi_n^\vartriangle x_k$, we have
$\{\|\widehat\psi_{k,n} - \psi_{k,n}^*\|^\vartriangle\}^2 = \{\|x_k
-\widehat\psi_{k,n}\|^\vartriangle\}^2 -
\{\|x_k-\psi_{k,n}^*\|^\vartriangle\}^2$ and
$\{\|x_k-\widehat\psi_{k,n}\|_n^\vartriangle\}^2 \leq\{\|x_k -
\psi_{k,n}^*\|_n^\vartriangle\}^2$. These two relations and Lemma
S.3 in the supplementary note imply that
$\|\widehat\psi_{k,n}-\psi^*_{k,n}\|^\vartriangle= \RMo_P(1)$,
which in turn by the same lemma implies
$\|\widehat\psi_{k,n}-\psi^*_{k,n}\|^\vartriangle_n = \RMo_P(1)$. As a
consequence, $\|\widehat\psi_{k,n}-\psi^*_k\|^\vartriangle_n = \RMo_P(1)$,
which is exactly~(\ref{eqapproxI-0}). The proof is complete.\vspace*{-2pt}
\end{appendix}

\section*{Acknowledgements}\vspace*{-2pt}

G. Cheng supported by NSF Grant DMS-09-06497 and NSF CAREER Award DMS-1151692.
L. Zhou supported in part by NSF Grant DMS-09-07170.
J. Z. Huang supported in part by NSF Grants DMS-06-06580, DMS-09-07170,
NCI (CA57030), and Award Number KUS-CI-016-04, made by King Abdullah University
of Science and Technology (KAUST).\vspace*{-2pt}

\begin{supplement}
\stitle{Supplement to ``Efficient semiparametric estimation in
generalized partially
linear additive models for longitudinal/clustered data''\\}
\slink[doi]{10.3150/12-BEJ479SUPP} 
\sdatatype{.pdf}
\sfilename{BEJ479\_supp.pdf}
\sdescription{The supplementary file (Cheng, Zhou and Huang~\cite
{chz12}) includes\vadjust{\goodbreak} the properties of the least favorable directions and
the complete proofs of Theorems \ref{thmrate} and \ref{thmnormality}
together with some empirical processes results for the
clustered/longitudinal data. The results of a simulation study that
compares our method with that by Carroll \textit{et al.} \cite{cmmy09}
are also included.}
\end{supplement}


\printhistory


\begin{thebibliography}{27}

\bibitem{bkrw93}
\begin{bbook}[mr]
\bauthor{\bsnm{Bickel},~\bfnm{Peter~J.}\binits{P.J.}},
  \bauthor{\bsnm{Klaassen},~\bfnm{Chris A.~J.}\binits{C.A.J.}},
  \bauthor{\bsnm{Ritov},~\bfnm{Ya'acov}\binits{Y.}} \AND
  \bauthor{\bsnm{Wellner},~\bfnm{Jon~A.}\binits{J.A.}}
(\byear{1993}).
\btitle{Efficient and Adaptive Estimation for Semiparametric Models}.
\bseries{Johns Hopkins Series in the Mathematical Sciences}.
\blocation{Baltimore, MD}: \bpublisher{Johns Hopkins Univ. Press}.
\bid{mr={1245941}}
\bptok{imsref}%
\end{bbook}
\endbibitem

\bibitem{cmmy09}
\begin{barticle}[author]
\bauthor{\bsnm{Carroll},~\bfnm{Raymond~J.}\binits{R.J.}},
  \bauthor{\bsnm{Maity},~\bfnm{Arnab}\binits{A.}},
  \bauthor{\bsnm{Mammen},~\bfnm{Enno}\binits{E.}} \AND
  \bauthor{\bsnm{Yu},~\bfnm{Kyusang}\binits{K.}}
(\byear{2009}).
\btitle{Efficient semiparametric marginal estimation for partially linear
  additive model for longitudinal/clustered data}.
\bjournal{Statistics in BioSciences}
\bvolume{1}
\bpages{10--31}.
\bptok{imsref}%
\end{barticle}
\endbibitem

\bibitem{c88}
\begin{barticle}[mr]
\bauthor{\bsnm{Chen},~\bfnm{Hung}\binits{H.}}
(\byear{1988}).
\btitle{Convergence rates for parametric components in a partly linear model}.
\bjournal{Ann. Statist.}
\bvolume{16}
\bpages{136--146}.
\bid{doi={10.1214/aos/1176350695}, issn={0090-5364}, mr={0924861}}
\bptok{imsref}%
\end{barticle}
\endbibitem

\bibitem{cj06}
\begin{barticle}[mr]
\bauthor{\bsnm{Chen},~\bfnm{Kani}\binits{K.}} \AND
  \bauthor{\bsnm{Jin},~\bfnm{Zhezhen}\binits{Z.}}
(\byear{2006}).
\btitle{Partial linear regression models for clustered data}.
\bjournal{J. Amer. Statist. Assoc.}
\bvolume{101}
\bpages{195--204}.
\bid{doi={10.1198/016214505000000592}, issn={0162-1459}, mr={2268038}}
\bptok{imsref}%
\end{barticle}
\endbibitem

\bibitem{chz12}
\begin{barticle}[author]
\bauthor{\bsnm{Cheng},~\bfnm{Guang}\binits{G.}},
  \bauthor{\bsnm{Zhou},~\bfnm{Lan}\binits{L.}} \AND
  \bauthor{\bsnm{Huang},~\bfnm{Jianhua~Z.}\binits{J.Z.}}
(\byear{2014}).
\btitle{Supplement to ``Efficient semiparametric estimation in generalized partially
linear additive models for longitudinal/clustered data.''
DOI:\doiurl{10.3150/12-BEJ479SUPP}}.
\bptok{imsref}%
\end{barticle}
\endbibitem

\bibitem{d01}
\begin{bbook}[mr]
\bauthor{\bparticle{de} \bsnm{Boor},~\bfnm{Carl}\binits{C.}}
(\byear{2001}).
\btitle{A Practical Guide to Splines},
\bedition{revised} ed.
\bseries{Applied Mathematical Sciences}
\bvolume{27}.
\blocation{New York}: \bpublisher{Springer}.
\bid{mr={1900298}}
\bptok{imsref}%
\end{bbook}
\endbibitem

\bibitem{dhlz02}
\begin{bbook}[mr]
\bauthor{\bsnm{Diggle},~\bfnm{Peter~J.}\binits{P.J.}},
  \bauthor{\bsnm{Heagerty},~\bfnm{Patrick~J.}\binits{P.J.}},
  \bauthor{\bsnm{Liang},~\bfnm{Kung-Yee}\binits{K.Y.}} \AND
  \bauthor{\bsnm{Zeger},~\bfnm{Scott~L.}\binits{S.L.}}
(\byear{2002}).
\btitle{Analysis of Longitudinal Data},
\bedition{2nd} ed.
\bseries{Oxford Statistical Science Series}
\bvolume{25}.
\blocation{Oxford}: \bpublisher{Oxford Univ. Press}.
\bid{mr={2049007}}
\bptok{imsref}%
\end{bbook}
\endbibitem

\bibitem{hlg00}
\begin{bbook}[author]
\bauthor{\bsnm{H{\"a}rdle},~\bfnm{Wolfgang}\binits{W.}},
  \bauthor{\bsnm{Liang},~\bfnm{Hua}\binits{H.}} \AND
  \bauthor{\bsnm{Gao},~\bfnm{Jiti}\binits{J.}}
(\byear{2000}).
\btitle{Partially Linear Models}.
\blocation{New York}: \bpublisher{Springer}.
\bptok{imsref}%
\end{bbook}
\endbibitem

\bibitem{hfz05}
\begin{barticle}[mr]
\bauthor{\bsnm{He},~\bfnm{Xuming}\binits{X.}},
  \bauthor{\bsnm{Fung},~\bfnm{Wing~K.}\binits{W.K.}} \AND
  \bauthor{\bsnm{Zhu},~\bfnm{Zhongyi}\binits{Z.}}
(\byear{2005}).
\btitle{Robust estimation in generalized partial linear models for clustered
  data}.
\bjournal{J. Amer. Statist. Assoc.}
\bvolume{100}
\bpages{1176--1184}.
\bid{doi={10.1198/016214505000000277}, issn={0162-1459}, mr={2236433}}
\bptok{imsref}%
\end{barticle}
\endbibitem

\bibitem{hzf02}
\begin{barticle}[mr]
\bauthor{\bsnm{He},~\bfnm{Xuming}\binits{X.}},
  \bauthor{\bsnm{Zhu},~\bfnm{Zhong-Yi}\binits{Z.Y.}} \AND
  \bauthor{\bsnm{Fung},~\bfnm{Wing-Kam}\binits{W.K.}}
(\byear{2002}).
\btitle{Estimation in a semiparametric model for longitudinal data with
  unspecified dependence structure}.
\bjournal{Biometrika}
\bvolume{89}
\bpages{579--590}.
\bid{doi={10.1093/biomet/89.3.579}, issn={0006-3444}, mr={1929164}}
\bptok{imsref}%
\end{barticle}
\endbibitem

\bibitem{hwz02}
\begin{barticle}[mr]
\bauthor{\bsnm{Huang},~\bfnm{Jianhua~Z.}\binits{J.Z.}},
  \bauthor{\bsnm{Wu},~\bfnm{Colin~O.}\binits{C.O.}} \AND
  \bauthor{\bsnm{Zhou},~\bfnm{Lan}\binits{L.}}
(\byear{2002}).
\btitle{Varying-coefficient models and basis function approximations for the
  analysis of repeated measurements}.
\bjournal{Biometrika}
\bvolume{89}
\bpages{111--128}.
\bid{doi={10.1093/biomet/89.1.111}, issn={0006-3444}, mr={1888349}}
\bptok{imsref}%
\end{barticle}
\endbibitem

\bibitem{hzz07}
\begin{barticle}[mr]
\bauthor{\bsnm{Huang},~\bfnm{Jianhua~Z.}\binits{J.Z.}},
  \bauthor{\bsnm{Zhang},~\bfnm{Liangyue}\binits{L.}} \AND
  \bauthor{\bsnm{Zhou},~\bfnm{Lan}\binits{L.}}
(\byear{2007}).
\btitle{Efficient estimation in marginal partially linear models for
  longitudinal/clustered data using splines}.
\bjournal{Scand. J. Stat.}
\bvolume{34}
\bpages{451--477}.
\bid{doi={10.1111/j.1467-9469.2006.00550.x}, issn={0303-6898}, mr={2368793}}
\bptok{imsref}%
\end{barticle}
\endbibitem

\bibitem{k99}
\begin{bbook}[mr]
\bauthor{\bsnm{Kress},~\bfnm{Rainer}\binits{R.}}
(\byear{1999}).
\btitle{Linear Integral Equations},
\bedition{2nd} ed.
\bseries{Applied Mathematical Sciences}
\bvolume{82}.
\blocation{New York}: \bpublisher{Springer}.
\bid{doi={10.1007/978-1-4612-0559-3}, mr={1723850}}
\bptok{imsref}%
\end{bbook}
\endbibitem

\bibitem{czp10}
\begin{barticle}[mr]
\bauthor{\bsnm{Leng},~\bfnm{Chenlei}\binits{C.}},
  \bauthor{\bsnm{Zhang},~\bfnm{Weiping}\binits{W.}} \AND
  \bauthor{\bsnm{Pan},~\bfnm{Jianxin}\binits{J.}}
(\byear{2010}).
\btitle{Semiparametric mean-covariance regression analysis for longitudinal
  data}.
\bjournal{J. Amer. Statist. Assoc.}
\bvolume{105}
\bpages{181--193}.
\bnote{With supplementary material available online}.
\bid{doi={10.1198/jasa.2009.tm08485}, issn={0162-1459}, mr={2656048}}
\bptok{imsref}%
\end{barticle}
\endbibitem

\bibitem{lz86}
\begin{barticle}[mr]
\bauthor{\bsnm{Liang},~\bfnm{Kung~Yee}\binits{K.Y.}} \AND
  \bauthor{\bsnm{Zeger},~\bfnm{Scott~L.}\binits{S.L.}}
(\byear{1986}).
\btitle{Longitudinal data analysis using generalized linear models}.
\bjournal{Biometrika}
\bvolume{73}
\bpages{13--22}.
\bid{doi={10.1093/biomet/73.1.13}, issn={0006-3444}, mr={0836430}}
\bptok{imsref}%
\end{barticle}
\endbibitem

\bibitem{lc01a}
\begin{barticle}[mr]
\bauthor{\bsnm{Lin},~\bfnm{Xihong}\binits{X.}} \AND
  \bauthor{\bsnm{Carroll},~\bfnm{Raymond~J.}\binits{R.J.}}
(\byear{2001}).
\btitle{Semiparametric regression for clustered data}.
\bjournal{Biometrika}
\bvolume{88}
\bpages{1179--1185}.
\bid{doi={10.1093/biomet/88.4.1179}, issn={0006-3444}, mr={1872228}}
\bptok{imsref}%
\end{barticle}
\endbibitem

\bibitem{lc01b}
\begin{barticle}[mr]
\bauthor{\bsnm{Lin},~\bfnm{Xihong}\binits{X.}} \AND
  \bauthor{\bsnm{Carroll},~\bfnm{Raymond~J.}\binits{R.J.}}
(\byear{2001}).
\btitle{Semiparametric regression for clustered data using generalized
  estimating equations}.
\bjournal{J. Amer. Statist. Assoc.}
\bvolume{96}
\bpages{1045--1056}.
\bid{doi={10.1198/016214501753208708}, issn={0162-1459}, mr={1947252}}
\bptok{imsref}%
\end{barticle}
\endbibitem

\bibitem{s92}
\begin{barticle}[mr]
\bauthor{\bsnm{Sasieni},~\bfnm{Peter}\binits{P.}}
(\byear{1992}).
\btitle{Nonorthogonal projections and their application to calculating the
  information in a partly linear {C}ox model}.
\bjournal{Scand. J. Stat.}
\bvolume{19}
\bpages{215--233}.
\bid{issn={0303-6898}, mr={1183198}}
\bptok{imsref}%
\end{barticle}
\endbibitem

\bibitem{s81}
\begin{bbook}[mr]
\bauthor{\bsnm{Schumaker},~\bfnm{Larry~L.}\binits{L.L.}}
(\byear{1981}).
\btitle{Spline Functions: Basic Theory}.
\blocation{New York}: \bpublisher{Wiley}.
\bid{mr={0606200}}
\bptok{imsref}%
\end{bbook}
\endbibitem

\bibitem{ss94}
\begin{barticle}[mr]
\bauthor{\bsnm{Severini},~\bfnm{Thomas~A.}\binits{T.A.}} \AND
  \bauthor{\bsnm{Staniswalis},~\bfnm{Joan~G.}\binits{J.G.}}
(\byear{1994}).
\btitle{Quasi-likelihood estimation in semiparametric models}.
\bjournal{J.~Amer. Statist. Assoc.}
\bvolume{89}
\bpages{501--511}.
\bid{issn={0162-1459}, mr={1294076}}
\bptok{imsref}%
\end{barticle}
\endbibitem

\bibitem{s88}
\begin{barticle}[mr]
\bauthor{\bsnm{Speckman},~\bfnm{Paul}\binits{P.}}
(\byear{1988}).
\btitle{Kernel smoothing in partial linear models}.
\bjournal{J. R. Stat. Soc. Ser. B Stat. Methodol.}
\bvolume{50}
\bpages{413--436}.
\bid{issn={0035-9246}, mr={0970977}}
\bptok{imsref}%
\end{barticle}
\endbibitem

\bibitem{s94}
\begin{barticle}[mr]
\bauthor{\bsnm{Stone},~\bfnm{Charles~J.}\binits{C.J.}}
(\byear{1994}).
\btitle{The use of polynomial splines and their tensor products in multivariate
  function estimation}.
\bjournal{Ann. Statist.}
\bvolume{22}
\bpages{118--171}.
\bptok{imsref}%
\end{barticle}
\endbibitem

\bibitem{vg00}
\begin{bbook}[author]
\bauthor{\bparticle{van~de} \bsnm{Geer},~\bfnm{Sara}\binits{S.}}
(\byear{2000}).
\btitle{Empirical Processes in M-Estimation}.
\blocation{Cambridge}: \bpublisher{Cambridge Univ. Press}.
\bptok{imsref}%
\end{bbook}
\endbibitem

\bibitem{w03}
\begin{barticle}[mr]
\bauthor{\bsnm{Wang},~\bfnm{Naisyin}\binits{N.}}
(\byear{2003}).
\btitle{Marginal nonparametric kernel regression accounting for within-subject
  correlation}.
\bjournal{Biometrika}
\bvolume{90}
\bpages{43--52}.
\bid{doi={10.1093/biomet/90.1.43}, issn={0006-3444}, mr={1966549}}
\bptok{imsref}%
\end{barticle}
\endbibitem

\bibitem{wcl05}
\begin{barticle}[mr]
\bauthor{\bsnm{Wang},~\bfnm{Naisyin}\binits{N.}},
  \bauthor{\bsnm{Carroll},~\bfnm{Raymond~J.}\binits{R.J.}} \AND
  \bauthor{\bsnm{Lin},~\bfnm{Xihong}\binits{X.}}
(\byear{2005}).
\btitle{Efficient semiparametric marginal estimation for longitudinal/clustered
  data}.
\bjournal{J. Amer. Statist. Assoc.}
\bvolume{100}
\bpages{147--157}.
\bid{doi={10.1198/016214504000000629}, issn={0162-1459}, mr={2156825}}
\bptok{imsref}%
\end{barticle}
\endbibitem

\bibitem{ZeDi94}
\begin{barticle}[auto:STB|2013/12/09|07:59:19]
\bauthor{\bsnm{Zeger},~\bfnm{S.~L.}\binits{S.L.}} \AND
  \bauthor{\bsnm{Diggle},~\bfnm{P.~J.}\binits{P.J.}}
(\byear{1994}).
\btitle{Semiparametric models for longitudinal data with application to CD4
  cell numbers in HIV seroconverters}.
\bjournal{Biometrics}
\bvolume{50}
\bpages{689--699}.
\bptok{imsref}%
\end{barticle}
\endbibitem

\bibitem{z04}
\begin{bmisc}[author]
\bauthor{\bsnm{Zhang},~\bfnm{L.}\binits{L.}}
(\byear{2004}).
\bhowpublished{Efficient estimation in marginal partially linear models for
  longitudinal/clustered data using splines.
Ph.D. thesis, Univ. Pennsylvania}.
\bptok{imsref}%
\end{bmisc}
\endbibitem

\end{thebibliography}
\end{document}